# EXACT AND APPROXIMATE RESULTS FOR DEPOSITION AND ANNIHILATION PROCESSES ON GRAPHS

By Mathew D. Penrose[1] and Aidan Sudbury

*University of Bath and Monash University*

We consider random sequential adsorption processes where the initially empty sites of a graph are irreversibly occupied, in random order, either by monomers which block neighboring sites, or by dimers. We also consider a process where initially occupied sites annihilate their neighbors at random times.

We verify that these processes are well defined on infinite graphs, and derive forward equations governing joint vacancy/occupation probabilities. Using these, we derive exact formulae for occupation probabilities and pair correlations in Bethe lattices. For the blocking and annihilation processes we also prove positive correlations between sites an even distance apart, and for blocking we derive rigorous lower bounds for the site occupation probability in lattices, including a lower bound of $1/3$ for $Z^2$. We also give normal approximation results for the number of occupied sites in a large finite graph.

**1. Introduction.** This article is concerned with certain Markov processes with state space $\{0,1\}^V$, where $V$ is the vertex set of a graph, in which the evolution at each site of $V$ is monotonic in time and the stochastic rate of evolution at a site is either a constant or zero, depending on the local configuration. Unusually for interacting particle systems, the focus is not on equilibrium measures (owing to the irreversibility) but rather on exact computation or estimation of probabilities, in which more progress than usual can be made because of the simple form of the interactions.

Two of the processes which we consider are versions of *random sequential adsorption* (RSA), in which the graph represents an initially "vacant" surface onto which particles are deposited successively at random, subject to

Received February 2004; revised April 2004.
[1]Supported in part by the Isaac Newton Institute for Mathematical Sciences, Cambridge.

an excluded volume effect whereby a deposited particle prevents new particles from being deposited nearby. We consider *monomer RSA with nearest-neighbor exclusion/blocking* (or *blocking RSA* for short) and *dimer RSA*. We also consider two forms of *annihilation process* in which sites are initially alive and from time to time kill their neighbors.

Interest in RSA from the physical sciences community is considerable: see the 1993 survey of Evans [7], and more recently, the series of surveys in Volume 165 (2000) of the journal *Colloids and Surfaces A*. The principal motivation comes from the notion of particles being physically deposited onto a surface. Motivation is discussed further in [19], where qualitative results are given for a class of process which generalizes those considered here on the integer lattice. Continuum versions of RSA are also of interest (see, e.g., [18, 20]) but are not considered here. The annihilation process was introduced by O'Hely and Sudbury in [14], as a model for the thinning of seedlings, where a plant can grow only by eliminating its neighbors.

In the present paper we consider general graphs (not only lattices). On the qualitative side, we derive an existence result for these processes on infinite graphs, forward equations governing the rate of change and a normal approximation theorem for the sum of the states of the vertices of a large finite graph. On the quantitative side, we derive explicit formulae for the occupation probability of sites in regular infinite trees (Bethe lattices) and also for correlations between sites in such trees. Moreover, we show that in the case of blocking RSA, these formulae provide good rigorous bounds for occupation probabilities in lattices such as $\mathbb{Z}^d$, comparing favorably to previously obtained rigorous bounds. The proof of these results for blocking RSA uses a result on positive correlations between sites an even distance apart, which is of interest in its own right.

A notable feature of these processes is that they can be interpreted in a totally discrete manner. For finite graphs, they can be generated via a uniform probability measure on the possible orderings in which a finite collection of "events" associated with the vertices or edges of the graph can occur. For infinite graphs, they may be obtained by taking a limit of finite graphs. Moreover, blocking RSA and the annihilation process on graphs provide a method of generating random independent sets of vertices, providing an alternative to the "hard-core" model which has received recent attention from discrete mathematicians [2, 9, 22].

## 2. Definitions and main results.

2.1. *Definitions for finite graphs.* Let $V$ be the vertex set, and $E$ the edge set, of a graph (we consider only graphs with no loops or multiple edges). We start by defining our models in the case where $V$ is finite. Each of them is described as a continuous-time stochastic process $(X_t^V, t \geq 0)$ taking values in $\{0,1\}^V$, with $X_t^V(v) = 1$ having the interpretation that site $v \in V$ is occupied at time $t$, and $X_t^V(v) = 0$ having the interpretation that $v$ is vacant at time $t$. For each $v$, the evolution of $(X_t^V(v), t \geq 0)$ is monotonic, so that the value of $X_t^V(v)$ changes at most once during the entire evolution of the process.

*Blocking RSA.* Let the random variables $(\tau_v, v \in V)$ be independent and uniformly distributed on $(0,1)$, with $\tau_v$ representing the arrival time of a particle at site $v$. All sites are initially vacant, and site $v$ becomes occupied at time $\tau_v$ if and only if none of its neighbors is already occupied. If site $v$ has one or more occupied neighbors at time $\tau_v$, the particle arriving at $v$ is rejected (blocked) and there is no change to the state of the system at time $\tau_v$ (this is the nearest-neighbor blocking/exclusion effect). Once occupied, a site remains permanently occupied.

*Dimer RSA.* Suppose $(\tau_e, e \in E)$ are independent and uniform on $(0,1)$. Here $\tau_e$ represents the arrival time of a dimer at edge $e$. All sites are initially vacant.

Let $e = \{u,v\} \in E$ be an edge of the graph, with endpoints $u, v$. If sites $u$ and $v$ are both vacant at time $\tau_e$, then the dimer arriving at edge $e$ at that time is accepted and sites $u, v$ both become occupied at that instant. If sites $u, v$ are not both vacant at time $\tau_e$, then the dimer arriving at $e$ is rejected, and nothing changes at time $\tau_e$. Once occupied, a site remains permanently occupied.

*The annihilation process.* Suppose $(\tau_e, e \in E)$ are independent and uniform on $(0,1)$. Here $\tau_e$ represents the time of occurrence of an "attack" along the edge $e$.

In this case each site is initially occupied. For each edge $e = \{u,v\}$, if sites $u, v$ are both occupied immediately prior to time $\tau_e$, then at time $\tau_e$ either site $u$ or site $v$ (each with probability $1/2$) becomes vacant (i.e., one site "attacks" the other). If $u, v$ are not both occupied at time $\tau_e$, then there is no change in the state of the system at that time. Once it becomes vacant, a site remains vacant permanently.

*The multiple annihilation process (MAP).* All sites are initially occupied, and $(\tau_v, v \in V)$ are independent and uniform on $(0,1)$. At time $\tau_v$, site $v$ (if still occupied) annihilates all of its neighbors and makes them all vacant, while remaining occupied itself. If $v$ is already vacant at time $\tau_v$ there is no change to the state of the system at time $\tau_v$.

Although defined in terms of annihilations, the MAP is most closely related to the blocking RSA process. To demonstrate this relationship and aid understanding of the MAP, we introduce a further process $(Z_t^V, t \geq 0)$ with state-space $\{0, 1, 2\}^V$, which we call the *undecided site process* (or USP for short). Sites in states 0, 1, 2 are interpreted as vacant, occupied and undecided, respectively. In the USP, all sites are initially undecided (i.e., in state 2). If site $v$ is undecided immediately prior to time $\tau_v$, then at that instant site $v$ becomes occupied and all neighboring sites become (or remain) vacant. Otherwise, there is no change at time $\tau_v$.

With the undecided site process $Z_t^V$ defined in this way, for $v \in V$ define coupled $\{0,1\}^V$-valued processes $(X_t^V, t \geq 0)$ and $(Y_t^V, t \geq 0)$ by

(1)  $\qquad X_t^V(v) = Y_t^V(v) = Z_t^V(v) \qquad \text{if } Z_t^V(v) \in \{0, 1\},$

(2)  $\qquad X_t^V(v) = 1 \text{ and } Y_t^V(v) = 0 \qquad \text{if } Z_t^V(v) = 2.$

In other words, undecided sites count as occupied for the process $(X_t^V, t \geq 0)$ and as vacant for the process $(Y_t^V, t \geq 0)$. It is not hard to see that with this construction, the process $(X_t^V, t \geq 0)$ is a realization of the MAP and $(Y_t^V, t \geq 0)$ is a realization of the blocking RSA process. At the terminal time 1, there are no undecided sites left and therefore the ultimate state of each site is the same for the coupled MAP and blocking RSA process, that is, $X_1^V(v) = Y_1^V(v)$ for all $v$. In particular, the distribution at time 1 of the MAP and the blocking RSA process are identical. In view of this similarity, we shall treat the MAP only briefly (in Section 4.4).

2.2. *Remarks.* For each of the blocking RSA, dimer RSA and annihilation processes, the evolution of $(X_t^V, t \geq 0)$ is determined by the order of occurrence of the "event times" $\tau_v$ or $\tau_e$ (together with the directions of "attacks" in the case of the annihilation process). Since they are assumed independent and identically distributed (i.i.d.), the "events" (arrivals of monomers, of dimers or attacks depending on the model) are equally likely to occur in any possible order. If we took the event times to be i.i.d. with some other nonuniform continuous distribution function $F$, then this would still be the case. In other words, changing the distribution function $F$ amounts merely to a reparametrization of time. If the distribution $F$ were chosen to be exponential, then the time-evolution of $X_t^V = (X_t^V(v), v \in V)$ would be that of a time-homogeneous continuous-time Markov chain with state-space $\{0, 1\}^V$. However, in fact it turns out to be more convenient to assume instead that $F$ is the uniform distribution on the unit interval $(0, 1)$, in which case these processes run as non-time-homogeneous continuous-time Markov chains.

Note also that while we could allow events to occur repeatedly at the same location (e.g., monomer arrivals at the times of a Poisson process for each site), in each model only the first "event" at any given location can possibly have any effect (e.g., if a monomer arriving at a particular site is blocked, then any later monomers arriving there are also blocked). Hence, only the time of the first arrival is relevant, and as just mentioned, we can just as well (and do) take this to be uniformly rather than exponentially distributed.

Owing to the monotonicity of the evolution, each process ultimately reaches a *jammed* state by time 1. For blocking (resp. dimer) RSA, jamming means that no vertex (edge) is any longer available to accept incoming particles. For annihilation, jamming means there is no pair of adjacent occupied sites. Both blocking RSA and annihilation produce a random independent set (in the graph-theoretic sense) of occupied sites at time 1; for blocking RSA the set of occupied sites is independent at all times. Of interest are the joint probability distributions of the states of sites at jamming, and also the states of sites at intermediate times prior to jamming.

2.3. *Extension to infinite graphs.* Let us now consider our three main processes on a *countably infinite* graph. While it is not hard to envisage that the preceding definitions of our monomer, dimer and annihilation processes should carry through to the case where $V$ is countably infinite, we need to check that these processes remain well defined. Evans ([7], page 1302) appeals to general results from Liggett [11] to do this, but we prefer a more direct approach, since the results in [11] are stated for interacting particle systems on $\mathbb{Z}^d$ rather than general graphs, and for time-homogeneous Markov processes rather than inhomogeneous ones, and in any case the direct approach is simpler here.

A sufficient condition for our processes to be well defined on an infinite graph is obtained by assuming that $V$ has bounded degree. In this case, the following result shows that we can define the blocking, dimer or annihilation process on $V$ as a limit of processes on finite subgraphs.

For any sequence of sets $(V_n)_{n \geq 1}$, we define the limit $\liminf(V_n) := \bigcup_{n=1}^{\infty} \bigcap_{m=n}^{\infty} V_m$.

PROPOSITION 1. *Suppose $(V, E)$ is a countably infinite graph with vertex set $V$ and edge set $E$, with uniformly bounded vertex degrees.*

(i) *Suppose $\tau_v, v \in V$, are independent and uniformly distributed on $[0, 1]$. Then, except on an event of probability zero, it is the case that for all $v \in V$ and $t \in [0, 1]$ and for any sequence $(V_n)_{n \geq 1}$ of finite subsets of $V$ with $\liminf(V_n) = V$, for blocking RSA the limit*

$$X_t^V(v) := \lim_{n \to \infty} (X_t^{V_n}(v)) \tag{3}$$

*exists and does not depend on the choice of $(V_n)_{n \geq 1}$.*

(ii) *Suppose $\tau_e, e \in E$, are independent and uniformly distributed on $[0, 1]$. Then, except on an event of probability zero, it is the case that for all $v \in V$ and $t \in [0, 1]$ and for any sequence $(V_n)_{n \geq 1}$ of finite subsets of $V$ with $\liminf(V_n) = V$, for either dimer RSA or the annihilation process the limit (3) exists and does not depend on the choice of $(V_n)_{n \geq 1}$.*

We use (3) to *define* the $\{0,1\}^V$-valued process $(X_t, t \geq 0)$, for countably infinite $V$. The almost sure convergence in (3) implies convergence of finite-dimensional distributions in the sense that that for any finite set of pairs $(t_1, v_1), \ldots, (t_k, v_k)$ in $[0,1] \times V$, and any sequence $(V_n)_{n \geq 1}$ of finite subsets of $V$ with $\liminf(V_n) = V$, it is the case that

$$(X_{t_j}^{V_n}(v_j), 1 \leq j \leq k) \xrightarrow{\mathcal{D}} (X_{t_j}^V(v_j), 1 \leq j \leq k) \qquad \text{as } n \to \infty.$$

Both in the preceding result and throughout this paper, we are adopting the following notational conventions. Notation such as $V$, $V_n$ or $W$ strictly speaking denotes the vertex set of a graph, but often will also be used to denote the graph itself. We often have occasion to consider one graph embedded in another; when $V \subset W$ and $W$ is the vertex set of a graph, then (e.g., for RSA on $V$) we assume the graph structure on $V$ is that induced by the graph structure on $W$. In other words, two elements of $V$ are adjacent (in the graph $V$) if and only if they are adjacent in the graph $W$.

Another way to define the process $(X_t^V, t \geq 0)$ on infinite graphs is by a recursion which we illustrate only for the case of blocking RSA. If $(\tau_v, v \in V)$ has a local minimum at site $w$, then denote site $w$ an *accepted site* and denote all neighbors of $w$ as *rejected sites*. Remove all accepted and rejected sites from $V$ and repeat. Continue indefinitely to generate a set of accepted sites. For each site $v$, and $t \geq 0$, set $X_t^V(v) = 1$ if and only if $t \geq \tau_v$ and site $v$ is accepted; otherwise set $X_t^V(v) = 0$. It can be shown that if $V$ is countable with bounded degree, the process $X_t^V$ generated this way is identical to that given in Proposition 1. For a continuum blocking RSA model, a similar approach is used and described in more detail by Penrose ([18], page 160).

2.4. *Forward equations.* Given a vertex set $V$ of a graph, write $\text{dist}(v, w)$ for the graph distance between two elements of $V$. We introduce the notation

(4) $$\mathcal{N}_k(v, V) := \{w \in V : \text{dist}(v, w) \leq k\}, \qquad k > 0,$$

for neighborhoods of vertices in $V$, and write just $\mathcal{N}_k(v)$ for $\mathcal{N}_k(v, V)$ if the graph $V$ is clear in the context.

The key to quantitative results is the fact that we can obtain differential equations for the probability of vacancy at particular sites. We first give such a result for blocking RSA.

PROPOSITION 2. *Let $V$ be finite or countably infinite with bounded degree. Let $W$ be a finite subset of $V$. Then for $0 \leq t < 1$, for blocking RSA,*

(5) $$\frac{d}{dt} P[X_t^V(v) = 0, v \in W]$$
$$= -\sum_{v \in W} P[X_t^{V \setminus \{v\}}(w) = 0, w \in (W \cup \mathcal{N}_1(v, V)) \setminus \{v\}].$$

We turn to the corresponding differential equations for dimer RSA and the annihilation process. Consider first the special case of vertices of degree 1.

PROPOSITION 3. *Let $V$ be the finite or countably infinite vertex set of a connected graph with bounded degree and more than one edge. Let $W$ be a finite subset of $V$, and suppose each element of $W$ has degree 1. Then for $0 \leq t < 1$, for dimer RSA,*

$$\begin{aligned}(6)\quad &\frac{d}{dt}P[X_t^V(v) = 0, v \in W] \\ &= -\sum_{v \in W} P[X_t^{V \setminus \{v\}}(w) = 0, w \in (W \cup \mathcal{N}_1(v,V)) \setminus \{v\}]\end{aligned}$$

*and for the annihilation process,*

$$\begin{aligned}(7)\quad &\frac{d}{dt}P[X_t^V(v) = 1, v \in W] \\ &= -\frac{1}{2}\sum_{v \in W} P[X_t^{V \setminus \{v\}}(w) = 1, w \in (W \cup \mathcal{N}_1(v,V)) \setminus \{v\}].\end{aligned}$$

The next result shows that the general case can be reduced to the case just considered, of vertices of degree 1, by considering a modified graph $\widetilde{G}_W$, which is illustrated in Figure 1.

Given a graph $G = (V, E)$, and given a vertex $w \in V$ of degree $d(w)$, enumerate the vertices in $\mathcal{N}_1(w, V)$ as $n_1(w), \ldots, n_{d(w)}(w)$, and let $\widetilde{G}_{\{w\}} = (\widetilde{V}_{\{w\}}, \widetilde{E}_{\{w\}})$ be the graph that is obtained from $G$ by replacing the vertex $w$ with $d(w)$ vertices of degree 1, denoted $v_1^*, \ldots, v_{d(w)}^*$, with $v_i^*$ connected to $n_i(w)$ by an edge for $1 \leq i \leq d(w)$.

Given a finite subset $W$ of $V$, let $\widetilde{G}_W = (\widetilde{V}_W, \widetilde{E}_W)$ be the graph obtained from $V$ by successively replacing each vertex $w$ in $W$ by a collection of $d(w)$ vertices of degree 1 as described above (the order in which the vertices in $W$ are taken is immaterial). In other words, define $\widetilde{G}_W$ for $W \subseteq V$ of cardinality $k$, inductively for $k = 1, 2, 3, \ldots$, by

$$(8)\qquad \widetilde{G}_{W \cup \{w\}}(V) = \widetilde{(\widetilde{G}_W)}_{\{w\}}.$$

For finite $W \subseteq V$, let $W^*$ be the set of vertices in $\widetilde{V}_W$ which replace the original vertices in $W$, a total of $\sum_{w \in W} d(w)$ vertices. Both $\widetilde{G}_W$ and $W^*$ are illustrated in Figure 1.

PROPOSITION 4. *Let $G = (V, E)$ be finite or countably infinite with bounded degree. Let $W$ be a finite subset of $V$. Then for $0 \leq t \leq 1$, for dimer RSA,*

$$(9)\qquad P[X_t^V(v) = 0, v \in W] = P[X_t^{\widetilde{V}_W}(v) = 0, v \in W^*]$$

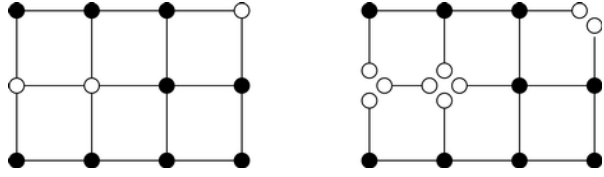

FIG. 1. *An example of the construction of $\widetilde{G}_W$. Here $W$ is the set of white vertices in the graph $G$ on the left, and $W^*$ is the set of white vertices in the graph $\widetilde{G}_W$ on the right.*

and for the annihilation process,

(10) $$P[X_t^V(v) = 1, v \in W] = P[X_t^{\widetilde{V}_W}(v) = 1, v \in W^*].$$

2.5. *Summary of exact results for regular trees.* Stochastic processes on trees have been the subject of much activity in the probability community in the last 10–20 years. See, for example, [12, 17]. One reason for looking at such processes is that they are sometimes more tractable than the equivalent processes on lattices such as $\mathbb{Z}^d$, while providing some insight into the latter processes.

In Section 4 we consider our RSA and annihilation processes on the $(k+1)$-regular tree (Bethe lattice), which we denote by $T_k$. This is a countably infinite connected graph in which each vertex has degree $k+1$ ($k \geq 1$ a fixed integer). We shall give exact computations of probabilities on $T_k$ which we summarize as follows.

For an arbitrary vertex $v$ of $T_k$, let $\alpha(t)$ denote the probability that vertex $v$ is still in its initial state (vacant in the case of blocking or dimer RSA on $T_k$, occupied in the case of annihilation on $T_k$) at time $t$. Let $\beta(t)$ be the probability that $v$ is still in its initial state at time $t$ for these processes on a graph obtained from $T_k$ by removing one branch from $v$ (in the case of blocking RSA) or by removing every branch except one from $v$ (in the case of dimer RSA or the annihilation process). We shall use the forward equations from the preceding section to derive the following formulae for $k \geq 2$ (only the second line needs modification for $k = 1$):

| Blocking RSA | Dimer RSA | Annihilation |
|---|---|---|
| $\beta'(t) = -\beta(t)^k$ | $\beta'(t) = -\beta(t)^k$ | $\beta'(t) = -\beta(t)^k/2$ |
| $\beta(t) = [1 + (k-1)t]^{-1/(k-1)}$ | $\beta(t) = [1 + (k-1)t]^{-1/(k-1)}$ | $\beta(t) = [1 + (k-1)t/2]^{-1/(k-1)}$ |
| $\alpha'(t) = -\beta(t)^{k+1}$ | $\alpha(t) = \beta(t)^{k+1}$ | $\alpha(t) = \beta(t)^{k+1}$ |
| $\alpha(t) = \frac{1}{2}(1 + \beta(t)^2)$ | | |

yielding the following numerical values for the occupation probability $P[X_t^{T_k}(v) = 1]$ at time $t = 1$:

| \multicolumn{4}{c}{**Occupation probabilities at time $t = 1$ on $T_k$**} |
| $k$ | **Blocking** | **Dimer** | **Annihilation** |
| --- | --- | --- | --- |
| 1  | 0.432 | 0.865 | 0.368 |
| 2  | 0.375 | 0.875 | 0.296 |
| 3  | 0.333 | 0.889 | 0.250 |
| 4  | 0.302 | 0.901 | 0.217 |
| 5  | 0.276 | 0.911 | 0.192 |
| 10 | 0.200 | 0.940 | 0.124 |
| 20 | 0.135 | 0.964 | 0.074 |

The formulae above for $\alpha(t)$ have previously been seen in the literature, not always in this form. In early work, Widom [23] gave an exact formula for the occupation density for blocking RSA in $T_2$ but later retracted this [24]. The correct formulae are given, at least for $t = 1$, by Evans ([5], page 2529) for dimer RSA, by Fan and Percus ([8], equation (4.8), see also [6]) for blocking RSA and by O'Hely and Sudbury [14] for the annihilation process. However, the first two of these derivations are quite involved, and probabilists may find our derivation more appealing.

We shall also give exact formulae for the covariances between the states of different sites of $T_k$, for each model. The formulae are given in Section 4, and in this summary we just give some of the resulting numerical values of correlations between two sites at time $t = 1$. Clearly these depend only on the graph distance between the two sites (and on $k$):

| | \multicolumn{6}{c}{**Correlations**} | | | | |
| | **Blocking** | | **Dimer** | | **Annihilation** | |
| **Distance** | $k = 3$ | $k = 5$ | $k = 3$ | $k = 5$ | $k = 3$ | $k = 5$ |
| --- | --- | --- | --- | --- | --- | --- |
| 1 | $-0.5000$ | $-0.3820$ | $-0.1250$ | $-0.0982$ | $-0.3333$ | $-0.2383$ |
| 2 | $0.1760$ | $0.1003$ | $0.0070$ | $0.0103$ | $0.0758$ | $0.0407$ |
| 3 | $-0.0474$ | $-0.0200$ | $0.0085$ | $0.0007$ | $-0.0130$ | $-0.0054$ |
| 4 | $0.0103$ | $0.0032$ | $-0.0041$ | $-0.0005$ | $0.0018$ | $0.0006$ |
| 5 | $-0.0019$ | $-0.0004$ | $0.0012$ | $0.0001$ | $-0.0002$ | $-0.0001$ |

See [16] for previous work on time-evolution of correlations for RSA on $T_1$. We believe that our formulae for correlations in $T_k, k \geq 2$, are new.

We note that the blocking and annihilation processes, but not the dimer process, show a pattern with sites at odd distances being negatively correlated and those at even distances being positively correlated. Next we note that correlations fall off extremely fast, since they are truncations of the exponential series (see Section 4). Also, correlations for the annihilation process are less than those for the blocking process. This is because for one site to influence another, not only do the "events" along a path between the two sites have to be in the correct temporal order, in the case of the annihilation process they have to be in the correct direction as well. Finally, the correlations for $k = 5$ are smaller than those for $k = 3$.

2.6. *Bounds for bipartite graphs.* A graph is *bipartite* if its vertex set can be partitioned into two independent (in the graph-theoretic sense) subsets. The next result says that for either RSA or the annihilation process, on a bipartite graph, the states of sites an even graph distance apart are positively correlated and the states of sites an odd graph distance apart are negatively correlated. This shows, among other things, that the pattern of positive and negative covariances for blocking RSA or the annihilation process, seen in the table from Section 2.5, extends to all $k$ and all distances (since any tree is bipartite).

THEOREM 5. *Suppose $V$ is the vertex set of a bipartite graph. For $j = 0, 1, \ldots, k$, suppose $W_j \subseteq V$, with $v_j \in W_j$ and $t_j \in [0,1]$, and define the event $E_j$ by*

$$E_j := \{X_{t_j}^{W_j}(v_j) = 0\} \qquad \text{if } \operatorname{dist}(v_j, v_0) \text{ is even},$$

$$E_j := \{X_{t_j}^{W_j}(v_j) = 1\} \qquad \text{if } \operatorname{dist}(v_j, v_0) \text{ is odd}.$$

*Then for either blocking RSA or the annihilation process,*

(11) $$P\left[\bigcap_{j=0}^{k} E_j\right] \geq P[E_0] P\left[\bigcap_{j=1}^{k} E_j\right].$$

In Theorem 5, it is to be understood that the processes $(X_t^{W_j}, t \geq 0)$ are defined for all $j$ in terms of the *same* set of random arrival times $(\tau_v, v \in V)$ (in the case of RSA) and the same set of random arrival times $\tau_e$ and random directions of attacks, along each edge (in the case of the annihilation process). The main case of interest in Theorem 5 is when $W_j = V$ for each $j$, but we shall also use the result in the more general form given.

Let us say $V$ is a *bipartite lattice* if (i) $V$ is the vertex set of a bipartite graph and (ii) for any $v_1, v_2, w_1, w_2 \in V$ with $v_1$ adjacent to $v_2$ and $w_1$ adjacent to $w_2$, there exists a graph isomorphism from $V$ to itself which sends $v_i$ to $w_i$, $i = 1, 2$. Examples include the integer lattice $\mathbb{Z}^d$ (with edges between each pair of vertices a unit Euclidean distance apart), and also the hexagonal lattice (but not the triangular lattice which is not bipartite). Various finite graphs are also bipartite lattices, including any lattice torus and also the discrete $n$-cube.

For a bipartite lattice, every vertex has the same degree. One of our main results says that for blocking RSA on such a lattice, the exact values given in the table for vacancy probabilities in the corresponding tree with the same vertex degrees, provide upper bounds for vacancy probabilities for the bipartite lattice.

THEOREM 6. *For blocking RSA on a bipartite lattice with vertices of degree $k+1$ where $k \geq 2$, for any vertex $v$ we have*

$$P[X_t^V(v) = 0] \leq \frac{1}{2}\left(1 + \frac{1}{(1+(k-1)t)^{2/(k-1)}}\right).$$

*In particular,*

$$P[X_1^V(v) = 0] \leq \tfrac{1}{2}(1 + k^{-2/(k-1)}) \qquad so \ P[X_1^V(v) = 1] \geq \tfrac{1}{2}(1 - k^{-2/(k-1)}).$$

In the special case of the square lattice $\mathbb{Z}^2$, we have $k = 3$ so that by Theorem 6, $P[X_t^V(v) = 0] \leq \frac{1+t}{1+2t}$ and the occupation probability at jamming satisfies $P[X_1^V(v) = 1] \geq 1/3$. Simulations suggest a value of approximately 0.364 ([7], page 1292).

In the case of the hexagonal lattice, we have $k = 2$ so that by Theorem 6, $P[X_1^V(v) = 1] \geq 3/8$, which compares favorably to the lower bound of 0.3596 provided (with a good deal of effort) by Caser and Hilhorst in [3]. Simulations indicate a value of 0.379 [6, 7, 23].

For the integer lattices $\mathbb{Z}^3$ and $\mathbb{Z}^4$ we have $k = 5$ and $k = 7$, respectively, so that Theorem 6 yields lower bounds for $P[X_1^V(v) = 1]$ of 0.2764 and 0.2386, respectively, compared to simulated values [7, 13] of 0.304 and 0.264, respectively.

REMARK. Our attempts to obtain good rigorous bounds for the annihilation process on bipartite lattices, using a similar argument to the proof of Theorem 6, have not been so successful. See the remarks at the end of the proof of Theorem 6.

2.7. *Results on normal approximation.* The next result is a normal approximation result for the total number of occupied sites at time $t$ on some (large) finite set $V$ for one of our processes on a (finite or infinite) graph with vertex set $W$ satisfying $V \subseteq W$. We denote this quantity $S_t(V, W)$, that is, set

(12) $$S_t(V, W) := \sum_{v \in V} X_t^W(v).$$

Let $\Phi(\cdot)$ denote the standard normal distribution function, and let $\text{card}(V)$ denote the number of elements (cardinality) of a finite set $V$.

12   M. D. PENROSE AND A. SUDBURYTHEOREM 7. *Let $D \geq 2$ be a finite constant and let $\varepsilon > 0$. Then there exists a finite constant $C = C(D, \varepsilon)$ such that for any finite $V$, any finite or countable $W$ such that $V \subseteq W$, and any graph with vertex set $W$ and maximum degree at most $D$, and all $x \in \mathbb{R}$, we have*

$$(13) \quad \left| P\left[\left(\frac{S_t(V,W) - ES_t(V,W)}{\sqrt{\mathrm{Var}(S_t(V,W))}}\right) \leq x\right] - \Phi(x) \right|$$
$$\leq C \mathrm{Var}(S_t(V,W))^{-(3/2)} (\mathrm{card}(V))^{1+\varepsilon}$$

*for any of blocking RSA, dimer RSA or the annihilation process.*

Suppose we are given a sequence of pairs $(V_n, W_n)_{n \geq 1}$, where for each $n$, $W_n$ is the vertex set of a finite or countable graph and $V_n \subseteq W_n$, and $\mathrm{card}(V_n) \to \infty$ as $n \to \infty$. For example, $W_n$ could be the same as $V_n$ for each $n$, or all $W_n$ could be the same infinite vertex set $W$, containing all the $V_n$. Under these circumstances, Theorem 7 gives us a central limit theorem with error bound on the rate of convergence, provided the $W_n$ have uniformly bounded degrees and we have reasonably fast growth of $\mathrm{Var}(S_t(V_n, W_n))$.

It can be shown that $\mathrm{Var}(S_t(V_n, W_n)) = O(\mathrm{card}(V_n))$, and one might expect $\mathrm{Var}(S_t(V_n, W_n)) = \Theta(\mathrm{card}(V_n))$ which would be sufficient to guarantee a central limit theorem from (13) with a rate of convergence $\mathrm{card}(V_n)^{\varepsilon - 1/2}$. The next result guarantees that this is indeed the case for $t < 1$.

THEOREM 8. *Let $t \in (0, 1)$. Let $D < \infty$. Then for any $V \subset W$ with all degrees in $W$ bounded by $D$, for blocking RSA,*

$$(14) \quad \mathrm{Var}(S_t(V,W)) \geq t(1-t)^{D+1}(D+1)^{-1} \mathrm{card}(V).$$

*If there are no isolated vertices in $V$, then for either dimer RSA or the annihilation process,*

$$(15) \quad \mathrm{Var}(S_t(V,W)) \geq (1/2)t(1-t)^{2D-1}(2D-1)^{-1} \mathrm{card}(V).$$

In the case $t = 1$, certain counterexamples suggest that, at least for blocking RSA or the annihilation process, there is no very simple condition to ensure that $\mathrm{Var}(S_1(V,W))/\mathrm{card}(V)$ is bounded away from zero. These are illustrated by Figure 2 and described below.

For the first example, suppose $V_0$ is an arbitrary finite graph, and $V$ is obtained by starting with $V_0$ and adding a "twin" vertex for each vertex in $V_0$, with each added vertex connected by an edge only to its twin. Ultimately, for each pair of twinned sites exactly one of the twins will be occupied, so that the variance of $S_1(V,V)$ is zero, for either blocking RSA or the annihilation process, although $V$ can be arbitrarily large.

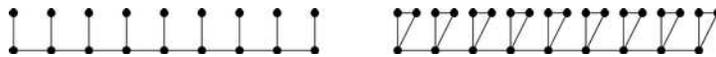

FIG. 2.  *Two graphs for which the ultimate number of particles deposited has zero variance.*

For a second example, suppose $V_0$ is an arbitrary finite graph, and $V$ is obtained by starting with $V_0$ and adding for each vertex $v \in V_0$ a pair of extra vertices $v', v''$, with $v'$ connected by edges only to $v''$ and to $v$, and $v''$ connected only to $v$ and to $v'$. In other words, a triangle is appended to each vertex in the original graph $V_0$, and again we conclude that the variance of $S_1(V,V)$ is zero for either blocking RSA or the annihilation process.

The first counterexample suggests that we may be in difficulties trying to find lower bounds for $\text{Var}(S_1(V,V))$ if all of $V$ lies near to the boundary (i.e., to the set of vertices of degree 1) while the second suggests that there can be problems if $V$ is not bipartite.

To give a sufficient condition for (14) to hold even for $t=1$, we introduce further terminology. Let us say that a site $v$ in a graph $W$ has *positive entropy* if (i) there exist at least two sites in $W$ that are adjacent to $v$ but not to each other and (ii) there exists an independent subset $B_v$ of $\mathcal{N}_3(v,W) \setminus \mathcal{N}_2(v,W)$ such that every vertex in $\mathcal{N}_2(v,W)$ has at least one neighbor in $B_v$. We shall later call $B_v$ a *blocking set* for $v$.

A sufficient condition for $v$ to have positive entropy is that $W$ be the vertex set of a bipartite graph and $v$ be a graph distance at least 3 from the boundary of $W$ [one can then take $B_v = \mathcal{N}_3(v,W) \setminus \mathcal{N}_2(v,W)$].

Denote the set of sites in $W$ with positive entropy by $W^+$.

THEOREM 9. *Let $D \in \mathbb{N}$. For any graph $W$ and any $V \subset W$ with all degrees in $W$ bounded by $D$, for blocking RSA,*

$$\text{Var}(S_1(V,W)) \geq \left( \frac{(1/2)^{D^3(1+D)+1}}{D^6(1+D)^2} \right) \text{card}(W^+ \cap V), \tag{16}$$

*and likewise there is a constant $\eta > 0$ depending only on $D$ such that for the annihilation process,*

$$\text{Var}(S_1(V,W)) \geq \eta \, \text{card}(W^+ \cap V). \tag{17}$$

It should be possible to find an analogous result to Theorem 9 for dimer RSA, but we have not done this.

We can easily deduce the following central limit theorem from the preceding results. Let $\mathcal{N}(0,\sigma^2)$ denote the normal distribution with mean zero and variance $\sigma^2$.

THEOREM 10. *Suppose $(W_n)_{n \geq 1}$ is a sequence of finite graphs with uniformly bounded degrees, and $V_n$ is a sequence of subgraphs of $W_n$ and with*

$\lim_{n \to \infty} \operatorname{card}(V_n) = \infty$. *Then for $0 < t < 1$,*

$$\text{(18)} \quad \frac{S_t(V_n, W_n) - ES_t(V_n, W_n)}{\sqrt{\operatorname{Var} S_t(V_n, W_n)}} \xrightarrow{\mathcal{D}} \mathcal{N}(0, 1)$$

*for any of blocking RSA, dimer RSA or the annihilation process. If, in addition, the sequence of pairs $(V_n, W_n)$ satisfies*

$$\liminf_{n \to \infty} (\operatorname{card}(V_n \cap W_n^+) / \operatorname{card}(V_n)) > 0,$$

*then* (18) *also holds for $t = 1$ in the case of either blocking RSA or the annihilation process.*

Theorem 10 has some overlap with Theorem 2.2 of [19]. However, Theorem 10 addresses processes on general graphs whereas [19] restricts attention to certain graphs embedded in $\mathbb{Z}^d$. Also, the proof of Theorem 10 provides a bound on the rate of convergence, unlike the corresponding proof in [19]. On the other hand, [19] is concerned with a more general class of processes than those considered here, and provides information on convergence of rescaled variances (for the graphs covered by [19]).

Finally, we can give more explicit central limit theorems in the one-dimensional case. Let $\mathbb{Z}$ be endowed with the usual nearest-neighbor graph structure (so that it is isomorphic to $T_1$ from Section 2.5). Let $L_n$ denote the set $\{1, 2, \ldots, n\}$, that is, a line of $n$ adjacent sites in $\mathbb{Z}$.

THEOREM 11. *Let $t \in (0, 1]$, and let $\sigma_{\mathrm{BL}}^2(t) := te^{-4t}$. For blocking RSA, as $n \to \infty$*

$$\text{(19)} \quad n^{-1} \operatorname{Var} S_t(L_n, \mathbb{Z}) \to \sigma_{\mathrm{BL}}^2(t),$$

$$\text{(20)} \quad n^{-1} \operatorname{Var} S_t(L_n, L_n) \to \sigma_{\mathrm{BL}}^2(t),$$

$$\text{(21)} \quad n^{-1}(S_t(L_n, \mathbb{Z}) - ES_t(L_n, Z)) \xrightarrow{\mathcal{D}} \mathcal{N}(0, \sigma_{\mathrm{BL}}^2(t))$$

*and*

$$\text{(22)} \quad n^{-1}(S_t(L_n, L_n) - ES_t(L_n, L_n)) \xrightarrow{\mathcal{D}} \mathcal{N}(0, \sigma_{\mathrm{BL}}^2(t)).$$

*Moreover, for dimer RSA and for the annihilation process the results* (19)–(22) *still hold with the limiting variance $\sigma_{\mathrm{BL}}^2$ replaced by $\sigma_{\mathrm{DI}}^2$ and $\sigma_{\mathrm{AN}}^2$, respectively, given by*

$$\sigma_{\mathrm{DI}}^2(t) := 4te^{-4t}, \qquad \sigma_{\mathrm{AN}}^2(t) := (1 + 2t)e^{-2t} - e^{-t}.$$

We note that using different methods, (19) was previously obtained, for blocking RSA only, by Pedersen and Hemmer [16], and the special case $t = 1$ of (20), for dimer RSA only, was obtained long ago by Page [15].

## 3. Proof of existence results and forward equations.

*Terminology.* By a *path of vertices* in a graph we mean a sequence of distinct vertices such that any two consecutive vertices in the sequence are connected by an edge. By a *path of edges*, we mean a sequence of distinct edges such that any two consecutive edges in the sequence have an endpoint in common.

PROOF OF PROPOSITION 1. (i) Suppose $W \subset W' \subset V$, with $W$ and $W'$ both finite. Then for $v \in W$, it is the case that $X_t^W(v) = X_t^{W'}(v)$ unless there is a path of vertices $v_0, v_1, \ldots, v_m$, starting at $v_0 = v$ and ending at some vertex $v_m \in W' \setminus W$, such that $\tau_{v_1} > \tau_{v_2} > \cdots > \tau_{v_m}$, since the only way that the effect of a vertex in $W' \setminus W$ can propagate to $v$ is along such a path in reverse order.

Therefore, for infinite $V$, a sufficient condition for the existence of the limit (3) for $v \in V$ is the existence of a finite random variable $R(v)$ such that there is no path $(v_0, v_1, \ldots, v_m)$ in $V$ with $v_0 = v$, $m > R(v)$, and with $\tau_{v_0} > \tau_{v_1} > \cdots > \tau_{v_m}$, since in that case, the value of $X_t^W(v)$ will be the same for all finite $W \subset V$ containing all vertices within graph distance $R(v)$ of $v$.

Let $v \in V$. Since degrees are assumed bounded by a constant which we denote $D$, the number of distinct self-avoiding paths starting at $v$ of length $m$ is bounded by $D^m$. For any such undirected path $(v_0, v_1, \ldots, v_m)$, we have

$$P[\tau_{v_0} > \tau_{v_1} > \cdots > \tau_{v_m}] = 1/(m+1)!$$

since any ordering of the times $\tau_{v_j}$, $0 \leq j \leq m$, is equally likely. Hence, the probability that there is a path of length $m$ starting from $v$ with associated arrival times in decreasing order is bounded by $D^m/(m+1)!$, and by the Borel–Cantelli lemma, the total number of paths starting from $v$ with associated arrival times in decreasing order is almost surely finite; we can take $R(v)$ to be the maximum of the lengths of such paths. This goes for any $v \in V$, and the existence of the limit in (3) for all $v, t$ follows.

(ii) Given $W \subset W' \subset V$, with $W$ and $W'$ finite, and given $v \in W$, it is the case that the value of $X_t^W(v) = X_t^{W'}(v)$ for dimer RSA, unless there exists a path of edges $e_1, e_2, \ldots, e_m$ that starts at $v$ and ends at some vertex in $W' \setminus W$, such that $\tau_{e_1} > \tau_{e_2} > \cdots > \tau_{e_m}$.

We assert that there exists almost surely a random variable $R'(v)$ such that there is no path of length greater than $R'(v)$, starting at $v$, with all edge arrival times along the path in decreasing order. The proof of this is the same as for the corresponding assertion in the proof of part (i).

Then $X_t^W(v)$ will be the same for all finite $W \subseteq V$ with $W$ containing all vertices a graph distance at most $R'(v)$ from $v$. Hence, (3) holds for the dimer process.

Next consider the annihilation process. With $v, W, W'$ as before, this time we have $X_t^W(v) = X_t^{W'}(v)$ unless there exists a path of edges $e_1, e_2, \ldots, e_m$ that starts at $v$ and ends at some vertex in $W' \setminus W$, such that $\tau_{e_1} > \tau_{e_2} > \cdots > \tau_{e_m}$ and in addition all "attacks" along the edges are in the direction of the path towards $v$. Hence with $R'(v)$ as defined for dimers, the annihilation process $X_t^W(v)$ will be the same for all finite $W \subseteq V$ with $W$ containing all vertices a graph distance at most $R'(v)$ from $v$. Hence, (3) holds for the annihilation process. $\square$

PROOF OF PROPOSITION 2. Suppose $V$ is finite or countably infinite with bounded degree, and $W$ is a finite subset of $V$. For $v \in W$ and $0 \leq t < t + h < 1$, we have

$$P[\{X_t^V(w) = 0, w \in W\} \cap \{X_{t+h}^V(v) = 1\}]$$
$$= P[\{\tau_v \in (t, t+h]\} \cap \{X_t^V(w) = 0, w \in (W \cup \mathcal{N}_1(v)) \setminus \{v\}\}] + o(h)$$
$$= hP[X_t^V(w) = 0, w \in (W \cup \mathcal{N}_1(v)) \setminus \{v\} | \tau_v \in (t, t+h]] + o(h)$$
$$= hP[X_t^{V \setminus \{v\}}(w) = 0, w \in (W \cup \mathcal{N}_1(v)) \setminus \{v\}] + o(h).$$

Summing over $v \in V$ we obtain

$$(23) \quad P\left[\{X_t^V(w) = 0, w \in W\} \cap \left(\bigcup_{v \in W} \{X_{t+h}(v) = 1\}\right)\right]$$
$$= h\left(\sum_{v \in W} P[X_t^{V \setminus \{v\}}(w) = 0, w \in (W \cup \mathcal{N}_1(v)) \setminus \{v\}]\right) + o(h).$$

Similarly, for $0 < h < t$,

$$P\left[\{X_{t-h}^V(w) = 0, w \in W\} \cap \left(\bigcup_{v \in W} \{X_t(v) = 1\}\right)\right]$$
$$(24) \quad = h\left(\sum_{v \in W} P[X_{t-h}^{V \setminus \{v\}}(w) = 0, w \in (W \cup \mathcal{N}_1(v)) \setminus \{v\}]\right) + o(h)$$
$$= h\left(\sum_{v \in W} P[X_t^{V \setminus \{v\}}(w) = 0, w \in (W \cup \mathcal{N}_1(v)) \setminus \{v\}]\right) + o(h),$$

where the last line follows from left-continuity in $t$ of $P[X_t^{V \setminus \{v\}}(w) = 0, w \in (W \cup \mathcal{N}_1(v)) \setminus \{v\}]$, which is easily seen. The result (5) follows from (23) and (24). $\square$

PROOF OF PROPOSITION 3. Let $V$ be the finite or countably infinite vertex set of a connected graph with bounded degree and more than one edge. Let $W$ be a finite subset of $V$, and suppose each element of $W$ has degree 1. Let $v \in W$. Then $\mathcal{N}_1(v)$ consists of a single vertex, denoted $u$, which is not itself in $W$ (otherwise, $G$ would either be disconnected or have only a single edge). Then for dimer RSA, writing $e$ for the edge $\{u,v\}$, we have

$$
\begin{aligned}
&P[\{X_t^V(w) = 0, w \in W\} \cap \{X_{t+h}^V(v) = 1\}] \\
(25) \quad &= P[\{\tau_e \in (t, t+h]\} \cap \{X_t^V(w) = 0, w \in (W \cup \mathcal{N}_1(v)) \setminus \{v\}\}] + o(h) \\
&= hP[X_t^{V \setminus \{v\}}(w) = 0, w \in (W \cup \mathcal{N}_1(v)) \setminus \{v\}] + o(h).
\end{aligned}
$$

Summing over $v \in W$, we obtain

$$
\begin{aligned}
&P\bigg[\{X_t^V(w) = 0, w \in W\} \cap \bigg(\bigcup_{v \in W}\{X_{t+h}(v) = 1\}\bigg)\bigg] \\
(26) \quad &= h \sum_{v \in W} P[X_t^{V \setminus \{v\}}(w) = 0, w \in (W \cup \mathcal{N}_1(v)) \setminus \{v\}] + o(h),
\end{aligned}
$$

and by arguing as for (24) we also have a similar expression to (26) for

$$
P\bigg[\{X_{t-h}^V(w) = 0, w \in W\} \cap \bigg(\bigcup_{v \in W}\{X_t(v) = 1\}\bigg)\bigg].
$$

Thus, we obtain (6) for dimer RSA.

For the annihilation process, given the event of an attack occurring along edge $\{u,v\}$ between times $t$ and $t+h$, this attack has conditional probability $1/2$ of being in the direction from $u$ to $v$, so that instead of (25) we obtain

$$
\begin{aligned}
&P[\{X_t^V(w) = 1, w \in W\} \cap \{X_{t+h}^V(v) = 0\}] \\
&= (h/2)P[X_t^{V \setminus \{v\}}(w) = 1, w \in (W \cup \mathcal{N}_1(v)) \setminus \{v\}] + o(h).
\end{aligned}
$$

The rest of the derivation of (7) for the annihilation process is the same as the derivation of (6) for dimer RSA. $\square$

PROOF OF PROPOSITION 4. Let $W$ be a finite subset of $V$. There is an obvious canonical bijection between the edges of $G$ and the edges of $\widetilde{G}_W$ (see Figure 1), and using this bijection, for either the dimer RSA or annihilation process there is a coupling between a realization of the process $(X_t^V, t \geq 0)$ and the process $(X_t^{\widetilde{V}_W}, t \geq 0)$ where we take the arrival time $\tau_e$ to be the same for edge $e$ of $G$ and for the corresponding edge of $\widetilde{G}_W$ (and in the annihilation model, we take the direction of attacks to be the same in both processes).

With this coupling, the state of each vertex $v \in V$ for the process $X_t^V$ is the same as the state of the corresponding vertex or vertices in $\widetilde{V}_W$ for the process $X_t^{\widetilde{V}_W}$, up to the instant when some site $v \in W$ first changes state, at which point one may see a difference between the coupled processes because only some of the sites in $\widetilde{V}_W$ change their state.

Thus, with this coupling of processes, the first instant at which some site in $W$ changes its state in the process $(X_t^V)$ is the same as the first instant at which some site in $W^*$ changes its state in the process $(X_t^{\widetilde{V}_W})$, and therefore, regardless of the coupling, the first instant at which some site in $W$ changes its state in the process $(X_t^V)$ has the same distribution as the first instant at which some site in $W^*$ changes its state in the process $(X_t^{\widetilde{V}_W})$. This is true for both dimer RSA and the annihilation process, so that (9) and (10) hold. $\square$

**4. Analysis of the $k$-tree.** This section is concerned with deriving the results for the case where $V$ is the vertex set of the $(k+1)$-regular tree $T_k$, which were summarized in Section 2.5.

It is clear that for each of the processes under consideration, for vertices $v, w$ of $T_k$ the covariance of $X_t^{T_k}(v)$ and $X_t^{T_k}(w)$ depends on $v, w$ only via $\mathrm{dist}(v, w)$. We shall use the notation

$$(27) \quad \mathrm{Cov}(m, t, T_k) := \mathrm{Cov}(X_t^{T_k}(v), X_t^{T_k}(w)) \qquad \text{if } \mathrm{dist}(v, w) = m.$$

Naturally $\mathrm{Cov}(m, t, T_k)$ will also depend on whether we are considering the blocking RSA, dimer, annihilation or multiple annihilation process.

Let $T_k^*$ be the rooted $(k+1)$-tree, that is, the graph with a single vertex (the root, denoted $v_0$) of degree $k$ and all others of degree $k+1$.

4.1. *Blocking RSA on $T_k$.*

THEOREM 12. *For blocking RSA on $T_k$, it is the case that for arbitrary $v \in T_k$,*

$$(28) \qquad P[X_t^{T_k}(v) = 0] = \tfrac{1}{2}(1 + (P[X_t^{T_k^*}(v_0) = 0])^2),$$

*with*

$$(29) \qquad P[X_t^{T_k^*}(v_0) = 0] = \begin{cases} (1 + (k-1)t)^{-1/(k-1)}, & k \geq 2, \\ e^{-t}, & k = 1. \end{cases}$$

*In particular,*

$$(30) \qquad P[X_1^{T_k}(v) = 0] = \begin{cases} \tfrac{1}{2}(1 + k^{-2/(k-1)}), & k \geq 2, \\ \tfrac{1}{2}(1 + e^{-2}), & k = 1. \end{cases}$$

PROOF. Fix $k$, and define $\beta(t)$ and $\alpha(t)$ (for an arbitrary vertex $v \in T_k$) by

$$\beta(t) = P[X_t^{T_k^*}(v_0) = 0], \qquad \alpha(t) := P[X_t^{T_k}(v) = 0] \tag{31}$$

(clearly the second probability does not depend on the choice of $v$). The subgraph of $T_k^*$ induced by removing $v_0$ consists of $k$ unconnected copies of $T_k^*$, so by Proposition 2,

$$\beta'(t) = -\beta(t)^k \tag{32}$$

with initial condition $\beta(0) = 1$, together giving us $\beta(t) = e^{-t}$ for $k = 1$, and

$$\frac{d}{dt}(\beta(t)^{1-k}) = k - 1 \quad \Longrightarrow \quad \beta(t)^{1-k} = 1 + (k-1)t, \qquad k \geq 2,$$

and (29) follows.

In a similar manner to (32), we obtain

$$\begin{aligned}\alpha'(t) &= -\beta(t)^{k+1} \\ &= \beta(t)\beta'(t) = \frac{d}{dt}\left(\frac{\beta(t)^2}{2}\right)\end{aligned} \tag{33}$$

with $\alpha(0) = \beta(0) = 1$, so that $\alpha(t) - 1 = (\beta(t)^2 - 1)/2$ and hence $\alpha(t) = (\beta(t)^2 + 1)/2$, that is, (28) holds. □

We now turn our attention to correlations between the states of two sites in $T_k$. Recall the definition (27).

THEOREM 13. *Let $k \geq 1$ be an integer, and consider blocking RSA on $T_k$. Set $\beta(t) := P[X_t^{T_k^*}(v_0) = 0]$. Then for $0 \leq t \leq 1$, and for integer $m \geq 1$,*

$$\operatorname{Cov}(m, t, T_k) = -\frac{\beta(t)^2}{2} \sum_{n=0}^{\infty} \frac{(2\log \beta(t))^{m+1+2n}}{(m+1+2n)!}. \tag{34}$$

REMARK. In the case $k = 1$, using (29) we see that this agrees with the formula of Pedersen and Hemmer [16] (our time parameter $t$ corresponds to $1 - e^{-t}$ in their parametrization of time; see our remarks in Section 2.2).

REMARK. In the case where $\operatorname{dist}(v, w) = 1$, our formula for the covariance in Theorem 13 simplifies to $-(1 - \alpha(t))^2$, with $\alpha(t)$ given in (31). This can be seen directly since for monomer RSA on any graph, $P[X_t^V(v) = X_t^V(w) = 1] = 0$ for adjacent sites $v, w$.

PROOF OF THEOREM 13. Again define $\alpha(t)$ and $\beta(t)$ by (31). Set
$$u(t) = 2\log(\beta(t)).$$

If $k \geq 2$, let $T_k^{**}$ denote an infinite tree in which a single vertex (the root, denoted $v_0$) has degree $k - 1$ and all the other vertices have degree $k + 1$. If $k = 1$, let $T_k^{**}$ consist of a single isolated vertex. Set
$$\gamma(t) := P[X_t^{T_k^{**}}(v_0) = 0].$$

Then by Proposition 2 and (32),

(35) $$\gamma'(t) = -\beta(t)^{k-1} = u'(t)/2.$$

Since $\gamma(0) = 1$ and $u(0) = 0$, we obtain

(36) $$\gamma(t) = 1 + u(t)/2.$$

For integer $m \geq 1$, let $T_{k,m}^{**}$ denote an infinite tree in which two vertices (denoted $v_0, v_1$) have degree $k$ and all other vertices have degree $k + 1$, and in which $\text{dist}(v_0, v_1) = m$. Set

$$\alpha_m(t) := P[X_t^{T_k}(v) = 0 \text{ and } X_t^{T_k}(w) = 0], \qquad v, w \in T_k, \ \text{dist}(v, w) = m,$$

$$\beta_m(t) := P[X_t^{T_k^*}(v_0) = 0 \text{ and } X_t^{T_k^*}(w) = 0], \qquad w \in T_k^*, \ \text{dist}(v_0, w) = m,$$

$$\gamma_m(t) := P[X_t^{T_{k,m}^{**}}(v_0) = 0 \text{ and } X_t^{T_{k,m}^{**}}(v_1) = 0].$$

Also set $\alpha_0(t) = \alpha(t)$, $\beta_0(t) = \beta(t)$, $\gamma_0(t) = \gamma(t)$. Then for $m \geq 1$, Proposition 2 followed by (35) yields

$$\gamma_m'(t) = -2\beta(t)^{k-1}\gamma_{m-1}(t) \quad \Longrightarrow \quad \frac{d\gamma_m}{du} = \gamma_{m-1},$$

with $\gamma_m(0) = 1$ and $u(0) = 0$, and with $\gamma_0(t) = 1 + u(t)/2$ by (36). This system of equations has solution

(37) $$\gamma_m = \frac{1}{2}\frac{u^{m+1}}{(m+1)!} + \sum_{n=0}^{m}\frac{u^n}{n!},$$

so that $\gamma_m(t) \to \beta(t)^2$ as $m \to \infty$, as would be expected.

Now consider $\beta_m(t)$. By Proposition 2 and (35), we have for $m \geq 1$ that

$$\beta_m'(t) = -\beta(t)^{k-1}\beta_{m-1}(t) - \beta(t)^k\gamma_{m-1}(t)$$

$$\Longrightarrow \quad \frac{d\beta_m}{du} = \frac{1}{2}(\beta_{m-1} + e^{u/2}\gamma_{m-1}),$$

with $\beta_m(0) = 1$ and $u(0) = 0$. This system of equations has solution

$$\beta_m = \frac{1}{2}\left(1 + \sum_{n=0}^{m}\frac{u^n}{n!}\right)e^{u/2}$$

and as $m \to \infty$ we have $\beta_m \to \frac{1}{2}(1 + \beta(t)^2)\beta(t) = \alpha(t)\beta(t)$, as expected.

Similarly for $\alpha_m(t)$, by Proposition 2 and (35) we have

$$(38) \quad \alpha'_m(t) = -2\beta(t)^k \beta_{m-1}(t) \implies \frac{d\alpha_m}{du} = e^{u/2}\beta_{m-1},$$

with $\alpha_m(0) = 1$ and $u(0) = 0$, giving us

$$\alpha_m = \frac{1+(-1)^m}{4} + \frac{e^u}{2}\left(1 + \sum_{n=0}^{\lfloor (m-1)/2 \rfloor} \frac{u^{m-1-2n}}{(m-1-2n)!}\right).$$

Hence,

$$\lim_{m \to \infty} \alpha_{2m} = \frac{1}{2} + \frac{e^u}{2}(1 + \sinh u)$$

$$= \frac{e^{2u} + 2e^u + 1}{4} = \frac{e^u}{2}(1 + \cosh u) = \lim_{m \to \infty} \alpha_{2m+1}$$

so that $\lim_{m \to \infty} \alpha_m(t) = \alpha(t)^2$, as one would expect. This implies the covariance between sites a distance $m$ apart is

$$\mathrm{Cov}(m, t, T_k) = \alpha_m(t) - \alpha(t)^2 = -\frac{1}{2}e^u \sum_{n=0}^{\infty} \frac{u^{m+1+2n}}{(m+1+2n)!},$$

which implies the asserted formula (34). □

4.2. *Dimer RSA on $T_k$.* In this section we derive the results for dimer RSA analogous to Theorems 12 and 13 from the preceding section, namely, (39) and (45) below for $k \geq 2$, and (40) and (46) for $k = 1$.

For dimer RSA on $T_k$, and for an arbitrary vertex $v$ of $T_k$, set

$$\alpha(t) = P[X_t^{T_k}(v) = 0], \qquad \beta(t) = P[X_t^{T_k^1}(v_0) = 0],$$

where $T_k^1$ has vertex $v_0$ of degree 1, all other vertices of degree $k + 1$. Then the graph $\widetilde{(T_k)}_{\{v\}}$ described at (8) has $k + 1$ components, each of them isomorphic to $T_k^1$. By Propositions 4 and 3,

$$\alpha(t) = \beta(t)^{k+1}, \qquad \beta'(t) = -\beta(t)^k,$$

and the solution to these [with $\beta(0) = 1$] is

$$(39) \quad \begin{aligned} \beta(t) &= (1 + (k-1)t)^{-1/(k-1)}, \\ \alpha(t) &= (1 + (k-1)t)^{-(k+1)/(k-1)}, \qquad k \geq 2, \end{aligned}$$

$$(40) \quad \beta(t) = e^{-t}, \qquad \alpha(t) = e^{-2t}, \qquad k = 1.$$

Set
$$\alpha_m(t) := P[X_t^{T_k}(v) = 0 \text{ and } X_t^{T_k}(w) = 0], \qquad v, w \in T_k, \text{ dist}(v, w) = m.$$

Define $T_{k,m}^1$ to be the tree with vertices $v_0, v_1$ of degree 1, with $\text{dist}(v_0, v_1) = m$, and all other vertices of degree $k + 1$. By Propositions 3 and 4 we have $\alpha_m(t) = \beta(t)^{2k} \beta_m(t)$, where
$$\beta_m(t) = P[X_t^{T_{k,m}^1}(v_0) = X_t^{T_{k,m}^1}(v_1) = 0].$$

Then $\beta_1(t) = 1 - t$, and by Propositions 4 and 3,
$$\beta_m'(t) = -2\beta(t)^{k-1}\beta_{m-1}(t). \tag{41}$$

To solve these, set $v(t) := 2\log\beta(t)$. Then $v'(t) = -2\beta(t)^{k-1}$ and
$$\frac{d\beta_m}{dv} = \frac{\beta_m'(t)}{v'(t)} = \beta_{m-1}, \qquad m \geq 2, \tag{42}$$

with $\beta_m(0) = 1$ and with
$$\beta_1 = 1 - t = \frac{k}{k-1} - \frac{e^{v(1-k)/2}}{k-1}, \qquad k \geq 2.$$

For any constants $a, b, c$, with $c \neq 0$, the solution to the system of equations (42) with $\beta_1(v) = a - be^{v/c}$ is
$$\beta_m = \frac{av^{m-1}}{(m-1)!} + \sum_{j=0}^{m-2} \frac{v^j}{j!} - bc^{m-1} \sum_{j=m-1}^{\infty} \frac{(v/c)^j}{j!}. \tag{43}$$

In the present case (for $k \geq 2$) we have $a = k/(k-1)$, $b = 1/(k-1)$ and $c = 2/(1-k)$ so that
$$\beta_m(t) = \frac{k(2\log\beta(t))^{m-1}}{(m-1)!(k-1)} + \left(\sum_{j=0}^{m-2} \frac{(2\log\beta(t))^j}{j!}\right)$$
$$+ \frac{1}{2}\left(\frac{2}{1-k}\right)^m \sum_{j=m-1}^{\infty} \frac{((1-k)\log\beta(t))^j}{j!}$$

and hence for $k \geq 2$,
$$\text{Cov}(m, t, T_k) \tag{44}$$
$$= \alpha_m(t) - \alpha(t)^2 = \beta(t)^{2k}(\beta_m(t) - \beta(t)^2)$$
$$= \beta(t)^{2k}\left(\frac{k(2\log\beta(t))^{m-1}}{(m-1)!(k-1)}\right. \tag{45}$$
$$\left. - \sum_{j=m-1}^{\infty}\left(\frac{(2\log\beta(t))^j}{j!} - \frac{1}{2}\left(\frac{2}{1-k}\right)^m \frac{((1-k)\log\beta(t))^j}{j!}\right)\right).$$

For $k=1$ we still have (41) with $\beta_1(t)=1-t$, but now we have $\beta(t)=e^{-t}$, so that $\beta(t)$ (for dimer RSA on $T_1$) equals $\gamma_{m-1}(t)$ (for monomer RSA on $T_1$, as given in Section 4.1), because they satisfy the same equations and agree for $m=1$. Thus

$$\beta_m(t) = \frac{1}{2}\frac{(-2t)^m}{m!} + \sum_{n=0}^{m-1}\frac{(-2t)^n}{n!},$$

and hence, since (44) still holds for $k=1$,

(46) $$\operatorname{Cov}(m,t,T_1) = -e^{-2t}\left\{\frac{1}{2}\frac{(-2t)^m}{m!} + \sum_{n=m+1}^{\infty}\frac{(-2t)^n}{n!}\right\}.$$

By (45) and (46), for all $k$ we have $\operatorname{Cov}(1,1,T_k) = -\alpha(1)^2$ which can also be seen directly since for adjacent sites $u,v$ we have $P[X_1^V(u) = X_1^V(v) = 0] = 0$ for dimer RSA.

4.3. *The annihilation process on $T_k$.* In this section we derive results for the annihilation process analogous to Theorems 12 and 13, namely, (47) and (49) below for $k \geq 2$, (48) and (50) for $k=1$. For the annihilation process, with $T_k^1$ as defined in the preceding section, let

$$\alpha(t) = P[X_t^{T_k}(v) = 1], \qquad \beta(t) = P[X_t^{T_k^1}(v_0) = 1].$$

Then Proposition 4 gives us $\alpha(t) = \beta(t)^{k+1}$, and Propositions 3 and 4 give us $\beta'(t) = -\beta(t)^k/2$. Also, $\beta(0) = 1$. Hence,

(47) $$\beta(t) = (1+(k-1)t/2)^{-1/(k-1)},$$
$$\alpha(t) = (1+(k-1)t/2)^{-(k+1)/(k-1)}, \qquad k \geq 2,$$

(48) $$\beta(t) = e^{-t/2}, \qquad \alpha(t) = e^{-t}, \qquad k=1.$$

Let $T_{k,m}^1$ be as described in the preceding section, and let

$$\alpha_m(t) := P[X_t^{T_k}(v) = X_t^{T_k}(w) = 1], \qquad \operatorname{dist}(v,w) = m,$$
$$\beta_m(t) := P[X_t^{T_{k,m}^1}(v_0) = X_t^{T_{k,m}^1}(v_1) = 1].$$

Then $\alpha_m(t) = \beta(t)^{2k}\beta_m(t)$, $\beta_1(t) = 1-t$, and by Propositions 3 and 4,

$$\beta'_m(t) = -\beta(t)^{k-1}\beta_{m-1}(t), \qquad m \geq 2.$$

Let $v(t) = 2\log\beta(t)$. Then $v'(t) = -\beta(t)^{k-1}$. Hence,

$$\frac{d\beta_m}{dv} = \beta_{m-1},$$

with

$$\beta_1(t) = 1 - t = \frac{k+1}{k-1} - \frac{2e^{(1-k)v/2}}{k-1}, \qquad k \geq 2.$$

Using the general solution (43) with $a = (k+1)/(k-1)$ and $b = 2/(k-1)$, and $c = 2/(1-k)$, for $k \geq 2$ we obtain

$$\beta_m(t) = \frac{(k+1)(2\log\beta(t))^{m-1}}{(m-1)!(k-1)} + \left(\sum_{j=0}^{m-2} \frac{(2\log\beta(t))^j}{j!}\right)$$
$$+ \left(\frac{2}{1-k}\right)^m \sum_{j=m-1}^{\infty} \frac{((1-k)\log\beta(t))^j}{j!},$$

so that for $k \geq 2$,

$$\text{Cov}(m,t,T_k)$$
$$= \alpha_m(t) - \alpha(t)^2 = \beta(t)^{2k}(\beta_m(t) - \beta(t)^2)$$
(49) $$= \beta(t)^{2k}\left(\frac{(k+1)(2\log\beta(t))^{m-1}}{(m-1)!(k-1)}\right.$$
$$\left. - \sum_{j=m-1}^{\infty}\left(\frac{(2\log\beta(t))^j}{j!} - \left(\frac{2}{1-k}\right)^m \frac{((1-k)\log\beta(t))^j}{j!}\right)\right).$$

In the case $k = 1$, we have $\beta'_m(t) = -\beta_{m-1}(t)$ with $\beta_1(t) = 1 - t$ and $\beta_m(0) = 1$, so that

$$\beta_m(t) = \sum_{j=0}^{m} \frac{(-t)^j}{j!}$$

and hence

(50) $$\text{Cov}(m,t,T_1) = \beta(t)^2(\beta_m(t) - \beta(t)^2) = -e^{-t} \sum_{j=m+1}^{\infty} \frac{(-t)^j}{j!}.$$

When $m = t = 1$, (49) and (50) yield $\text{Cov}(1,1,T_k) = -\alpha(1)^2$ for all $k \geq 1$. This can also be seen directly since for adjacent sites $u, v$ in any graph, we have $P[X_1^V(u) = X_1^V(v) = 1] = 0$ for the annihilation process.

4.4. *The multiple annihilation process on $T_k$, $k \geq 2$.* In this section we derive results for the multiple annihilation process (MAP) in $T_k$ analogous to Theorems 12 and 13 from Section 4.1, namely, (58) and (63). We consider only the case with $k \geq 2$.

We shall utilize the coupled construction of a realization of the MAP denoted $(X_t^V, t \geq 0)$ and a realization of the blocking RSA process denoted $(Y_t^V, t \geq 0)$ in terms of an undecided site process (USP) $(Z_t^V, t \geq 0)$ with state space $\{0,1,2\}^V$, as described in Section 2.1. For any vertices $v, w \in V$, the identities (1) and (2) give us

(51) $$P[Z_t^V(v) \in \{0,2\}] = P[Y_t^V(v) = 0],$$

(52) $P[\{Z_t^V(v) \in \{0,2\}\} \cap \{Z_t^V(w) \in \{0,2\}\}] = P[Y_t^V(v) = Y_t^V(w) = 0]$

and

(53) $$P[Z_t^V(v) \in \{1,2\}] = P[X_t^V(v) = 1],$$

(54) $P[\{Z_t^V(v) \in \{1,2\}\} \cap \{Z_t^V(w) \in \{1,2\}\}] = P[X_t^V(v) = X_t^V(w) = 1].$

With $T_k^*$ defined in Section 4.1, we use the following notation. For $i \in \{0,1,2\}$, define

$$\tilde{\beta}(i,t) := P[Z_t^{T_k^*}(v_0) = i], \qquad \tilde{\alpha}(i,t) := P[Z_t^{T_k}(v) = i].$$

The tilde indicates probabilities which refer to the USP. By (51),

(55) $$\tilde{\beta}(0,t) + \tilde{\beta}(2,t) = P[Y_t^{T_k^*}(v_0) = 0] = \beta(t),$$

(56) $$\tilde{\alpha}(0,t) + \tilde{\alpha}(2,t) = P[Y_t^{T_k}(v_0) = 0] = \alpha(t),$$

where $\beta(t)$ and $\alpha(t)$ are as defined in Section 4.1 (and refer to the blocking RSA process).

To remain "undecided" (i.e., in state 2) at time $t$, a site $v$ must have $\tau_v > t$ and must have all of its neighbors in state 0 or 2 at time $t$, so that by (55),

$$\tilde{\alpha}(2,t) = (1-t)(\tilde{\beta}(0,t) + \tilde{\beta}(2,t))^{k+1} = (1-t)\beta(t)^{k+1}.$$

Thus by (53) and (56), we have

(57) $$P[X_t^{T_k}(v) = 1] = \tilde{\alpha}(2,t) + \tilde{\alpha}(1,t) = (1-t)\beta(t)^{k+1} + 1 - \alpha(t).$$

Substituting the formulae for $\beta(t)$ and $\alpha(t)$ derived in Section 4.1 into (57), we obtain the following formula for the probability that a site is occupied for the MAP in $T_k, k \geq 2$:

(58) $$P[X_t^{T_k}(v) = 1] = (1-t)(1 + (k-1)t)^{(1+k)/(1-k)}$$
$$+ \tfrac{1}{2}(1 - (1 + (k-1)t)^{2/(1-k)}).$$

In order to calculate covariances, let $v, w$ be vertices of $T_k$ with $\text{dist}(v,w) = m$. For $i, j \in \{0,1,2\}$, define

$$\tilde{\alpha}_m(i,j,t) := P[Z_t^{T_k}(v) = i \text{ and } Z_t^{T_k}(w) = j].$$

For $v$ and $w$ to remain undecided (i.e., in state 2) at time $t$ requires $\tau_v > t$, $\tau_w > t$ and all neighbors in state 0 or 2 at time $t$. Hence by (51) and (52),

(59) $$\tilde{\alpha}_m(2,2,t) = (1-t)^2 \beta(t)^{2k} \gamma_{m-2}(t),$$

where $\gamma_m(t), m \geq 0$, is as defined in Section 4.1, and so is given by (37), and where we set $\gamma_{-1}(t) = 1$.

Next, we use elementary probability followed by (51) and (52) to obtain

$$\tilde{\alpha}_m(1,1,t) = 1 - P[Z_t^{T_k}(v) \in \{0,2\}] - P[Z_t^{T_k}(w) \in \{0,2\}]$$
(60)
$$+ P[\{Z_t^{T_k}(v) \in \{0,2\}\} \cap \{Z_t^{T_k}(w) \in \{0,2\}\}]$$
$$= 1 - 2\alpha(t) + \alpha_m(t),$$

with $\alpha(t)$ and $\alpha_m(t)$ given in Section 4.1.

We obtain $\tilde{\alpha}_m(1,2,t)$ indirectly, via a forward equation for $\tilde{\alpha}_m(1,1,t)$. As $h \downarrow 0$, we have

$$P[\{Z_{t+h}^{T_k}(v) = Z_{t+h}^{T_k}(w) = 1\} \setminus \{Z_t^{T_k}(v) = Z_t^{T_k}(w) = 1\}]$$
$$= 2P[\{Z_t^{T_k}(v) = 1\} \cap \{Z_t^{T_k}(w) = 2\} \cap \{t < \tau_w \leq t+h\}] + o(h)$$
$$= 2h(1-t)^{-1}P[\{Z_t^{T_k}(v) = 1\} \cap \{Z_t^{T_k}(w) = 2\}] + o(h),$$

the last line arising because the conditional distribution of $\tau_w$, given that event $\{Z_t^{T_k}(v) = 1\} \cap \{Z_t^{T_k}(w) = 2\}$ occurs, is uniform over the interval $(t,1]$. Thus,

$$\frac{d}{dt}\tilde{\alpha}_m(1,1,t) = 2\tilde{\alpha}_m(1,2,t)/(1-t).$$

Rearranging this identity and using (60) followed by (33) and (38), we obtain

(61)
$$\tilde{\alpha}_m(1,2,t) = \left(\frac{1-t}{2}\right)(-2\alpha'(t) + \alpha'_m(t))$$
$$= (1-t)(\beta(t)^{k+1} - \beta(t)^k \beta_{m-1}(t)),$$

with $\beta(t)$ and $\beta_m(t)$ given in Section 4.1.

Combining (54), (59), (60) and (61), we obtain

(62)
$$P[\{X_t^{T_k}(v) = 1\} \cap \{X_t^{T_k}(w) = 1\}]$$
$$= \tilde{\alpha}_m(1,1,t) + 2\tilde{\alpha}_m(1,2,t) + \tilde{\alpha}_m(2,2,t)$$
$$= 1 - 2\alpha(t) + \alpha_m(t) + 2(1-t)\beta(t)^k(\beta(t) - \beta_{m-1}(t))$$
$$+ (1-t)^2 \beta(t)^{2k}\gamma_{m-2}(t).$$

Let $\text{Cov}_{\text{MAP}}(m,t,T_k)$ and $\text{Cov}_{\text{RSA}}(m,t,T_k)$ be defined as at (27), with respect to the MAP and to blocking RSA, respectively. In particular, $\text{Cov}_{\text{RSA}}(m,t,T_k) = \alpha_m(t) - \alpha(t)^2$. Using (62) and (57), we obtain

$$\text{Cov}_{\text{MAP}}(m,t,T_k)$$
(63)
$$= P[\{X_t^{T_k}(v) = 1\} \cap \{X_t^{T_k}(w) = 1\}] - (P[X_t^{T_k}(v) = 1])^2$$
$$= \text{Cov}_{\text{RSA}}(m,t,T_k) + 2(1-t)\beta(t)^k(\alpha(t)\beta(t) - \beta_{m-1}(t))$$
$$+ (1-t)^2 \beta(t)^{2k}(\gamma_{m-2}(t) - \beta(t)^2).$$

In particular, $\text{Cov}_{\text{MAP}}(m, 1, T_k) = \text{Cov}_{\text{RSA}}(m, 1, T_k)$, which is as it should be.

## 5. Proof of bounds for bipartite graphs.

PROOF OF THEOREM 5. First we consider the blocking RSA process. Define the independent $U(0,1)$ variables $U_v, v \in V$, by

$$U_v = \begin{cases} \tau_v, & \text{if } \text{dist}(v, v_0) \text{ is even,} \\ 1 - \tau_v, & \text{if } \text{dist}(v, v_0) \text{ is odd.} \end{cases}$$

We assert that for $j = 0, 1, \ldots, k$, the indicator of the event $E_j$, defined by $E_j := \{X_{t_j}^{W_j}(v_j) = 0\}$ if $\text{dist}(v_0, v_j)$ is even, and by $E_j := \{X_{t_j}^{W_j}(v_j) = 1\}$ if $\text{dist}(v_0, v_j)$ is odd, is an increasing function of the variables $U_v, v \in V$. To see this choose any $v \in V$, let $0 \le u < u' \le 1$ and fix numbers $u_w \in [0, 1]$ for $w \in V \setminus \{v\}$. We assert that if

(64) $\quad \{U_v = u\} \cap \bigcap_{w \in V \setminus \{v\}} \{U_w = u_w\} \quad \text{implies } E_j \text{ occurs,}$

then

(65) $\quad \{U_v = u'\} \cap \bigcap_{w \in V \setminus \{v\}} \{U_w = u_w\} \quad \text{implies } E_j \text{ occurs.}$

To see this, first suppose $\text{dist}(v, v_0)$ is even. Consider the effect of increasing the value of $U_v$ (and hence of $\tau_v$) from $u$ to $u'$ while leaving the values of $U_w, w \in V \setminus \{v\}$, unchanged. This can change the value of $X_t^{W_j}(v)$, if at all, by changing it from 1 to 0. It can also directly affect sites $w$ neighboring $v$, if at all, by changing $X_t^{W_j}(w)$ from 0 to 1. These sites, in turn, can affect some of the sites $w$ at a distance 2 from $v$ by changing $X_t^{W_j}(w)$ from 1 to 0. In general, for any $w, t$, the effect of changing $U_v$ from $u$ to $u'$ on the value of $X_t^{W_j}(w)$, if any, can only be to decrease it if $\text{dist}(v, w)$ is even and can only be to increase it if $\text{dist}(v, w)$ is odd. The assertion that (64) implies (65) follows.

The argument is similar when $\text{dist}(v, v_0)$ is odd. In this case, increasing the value of $U_v$ is equivalent to decreasing the value of $\tau_v$, and this can only increase the value of $X_t^{W_j}(w)$ when $\text{dist}(v, w)$ is even and can only decrease the value of $X_t^{W_j}(w)$ when $\text{dist}(v, w)$ is odd. Since $\text{dist}(v, w)$ is even if and only if $\text{dist}(v_0, w)$ is odd, we may conclude once again that (64) implies (65). This justifies our initial assertion that the indicator of event $E_j$ is increasing in the variables $U_v$, for each $j$.

The result (11) for the blocking process now follows from the Harris–FKG inequality. The usual version of this result (see, e.g., [10]) gives positive correlations of events which are increasing in a family of $\{0, 1\}$-valued variables, but the proof in [10] is easily adapted to the present case with $[0,1]$-valued variables $U_j$.

Now consider the annihilation process. We first consider a modified version of this process, whereby for any pair of neighboring sites $u, v$ there is a random time $t_{uv}$ uniform on $[0,1]$ at which $u$ annihilates $v$ (if both sites were still occupied immediately before time $t_{uv}$). It is clear that only the event at the minimum of $t_{uv}, t_{vu}$ can be effective. Therefore, this modified annihilation process is equivalent to the one we originally defined, except that for each edge the time $\tau_e$ is distributed as the minimum of two independent uniform variables on $(0,1)$, rather than as a uniform variable on $(0,1)$. However, this, in turn, is equivalent to a time-change of the original annihilation process, as explained in Section 2.2, and therefore it suffices to prove the result for this modified annihilation process.

Let the uniform random variables $U_{u,v}$ be defined for each pair of neighboring sites $u, v$, by

$$U_{uv} = \begin{cases} t_{uv}, & \text{if } \operatorname{dist}(u, v_0) \text{ is even,} \\ 1 - t_{uv}, & \text{if } \operatorname{dist}(u, v_0) \text{ is odd.} \end{cases}$$

By the same argument as for blocking RSA, for each $j$ the indicator of event $E_j$ is an increasing function of the variables $U_{uv}$, and the proof of (11) for blocking RSA carries through to the annihilation process. □

PROOF OF THEOREM 6. Suppose $V$ is a bipartite lattice and each vertex has degree $k + 1$. Let $v_0$ be some arbitrary specified vertex of $V$, let the neighbors of $v_0$ be denoted $v_1, \ldots, v_{k+1}$ and let the neighbors of $v_1$ (other than $v_0$) be denoted $w_1, \ldots, w_k$. Let $V^*$ denote the graph $V \setminus \{v_0\}$, and let $V^{**}$ denote the graph $V \setminus \{v_0, v_1\}$.

For the blocking RSA process, define

$$\beta(t) := P[X_t^{V^*}(v_1) = 0],$$
$$\delta(t) := P[X_t^{V^{**}}(w_1) = X_t^{V^{**}}(w_2) = \cdots := X_t^{V^{**}}(w_k) = 0],$$
$$\eta(t) := P[X_t^{V^*}(v_1) = X_t^{V^*}(v_2) = \cdots = X_t^{V^*}(v_{k+1}) = 0].$$

If we let $\mathcal{E}$ denote the event that $w_1, \ldots, w_k$ are all unoccupied at time $t_1$, then $\delta(t)$ is the probability of event $\mathcal{E}$ for blocking RSA on $V \setminus \{v_1\}$ with the value of $\tau_{v_0}$ set to be (say) 2. Since $\mathbf{1}_{\mathcal{E}}$ is increasing in $\tau_{v_0}$, $\delta(t) \geq P[\mathcal{E}]$ for blocking RSA on $V \setminus \{v_1\}$ without any condition on $\tau_{v_0}$. Then by positive correlations between sites at an even graph distance (Theorem 5), we have

$$\delta(t) \geq P[\mathcal{E}] \geq \beta(t)^k$$

so that by the forward equation (Proposition 2),

(66) $$\beta'(t) = -\delta(t) \leq -\beta(t)^k,$$

with $\beta(0) = 1$. Hence,

(67)
$$\frac{d}{dt}\left(-\frac{1}{\beta(t)^{k-1}}\right) \leq 1 - k \implies -\frac{1}{\beta(t)^{k-1}} \leq (1-k)t - 1$$
$$\implies \beta(t) \leq \frac{1}{(1 + (k-1)t)^{1/(k-1)}}.$$

Next, we assert that

(68) $$\eta(t) \geq \delta(t)\beta(t).$$

To see this, note that $\eta(t)$ is the probability of the intersection of events $A \cap B$, where we set $A := \{X_t^{V^*}(v_{k+1}) = 0\}$ and

$$B := \{X_t^{V^*}(v_1) = X_t^{V^*}(v_2) = \cdots = X_t^{V^*}(v_k) = 0\},$$
$$B' := \{X_t^{V^*\setminus\{v_{k+1}\}}(v_1) = X_t^{V^*\setminus\{v_{k+1}\}}(v_2) = \cdots = X_t^{V^*\setminus\{v_{k+1}\}}(v_k) = 0\}.$$

Then $P[B'] = \delta(t)$.

It is to be understood that $A, B, B'$ are all defined in terms of the same set of random arrival times $\{\tau_v, v \in V\}$. In effect, $B'$ is the event that for RSA on $V^*$, if we increase the value of $\tau_{v_{k+1}}$ to 2 (say) but leave $\tau_v, v \in V^* \setminus \{v_{k+1}\}$, the same, then we have sites $v_1, \ldots, v_k$ vacant at time $t$. Then we have the event equality

$$A \cap B = A \cap B',$$

since if $A$ occurs, then increasing the value of $\tau_{v_{k+1}}$ to 2 does not have any effect on the acceptances or rejections up to time $t$ at sites other than $v_{k+1}$ (since the arrival at site $v_{k+1}$ was rejected even with the original value of $\tau_{v_{k+1}}$, at least if that value was at most $t$). Events $A$ and $B'$ are positively associated by Theorem 5, so (68) follows.

Setting $\alpha(t) = P[X_t^V(v_0) = 0]$, by the forward equation (Proposition 2), along with (68) and (66), we have

$$\alpha'(t) = -\eta(t) \leq -\delta(t)\beta(t) = \beta'(t)\beta(t).$$

Since $\alpha(0) = \beta(0) = 1$, we then have

$$\alpha(t) - 1 = \int_0^t \alpha'(s)\,ds \leq \int_0^t \beta(s)\beta'(s)\,ds = \tfrac{1}{2}(\beta(t)^2 - 1),$$

so that by (67),

$$\alpha(t) \leq \frac{1}{2}(1 + \beta(t)^2) \leq \frac{1}{2}\left(1 + \frac{1}{(1+(k-1)t)^{2/(k-1)}}\right). \qquad \square$$

REMARK. For the annihilation process, the argument up to a point is the same as for the blocking RSA process. We have

$$\alpha(t) \geq \beta(t)^{k+1}, \qquad \beta'(t) \leq -\frac{1}{2}\beta(t)^k \quad \Longrightarrow \quad \beta(t) \leq \left[\frac{2}{2+(k-1)t}\right]^{1/(k-1)}.$$

Unfortunately, the two inequalities appear to be in opposite directions. In fact, on square lattices simulations suggest that $\alpha(t) >$ the value on the tree (0.252 as opposed to 0.25). This may make this approach difficult.

**6. Proof of results relating to normal approximation.** To prove Theorem 7, we shall use the following theorem on normal approximation of sums of locally dependent random variables, which is a special case of Theorem 2.4 of Chen and Shao [4].

LEMMA 14 ([4]). *Suppose $V$ is a finite set and $\{\xi_v, v \in V\}$ is a collection of random variables, each with mean zero and finite variance. For $A \subseteq V$ set $\xi_A := \sum_{v \in A} X_v$, and suppose that $\mathrm{Var}(\xi_V) = 1$. Let $\kappa$ be a constant.*

*Suppose that for each $v \in V$ there exist subsets $A_v \subseteq B_v \subseteq C_v$ of $V$ such that $\xi_v$ is independent of $\xi_{V \setminus A_v}$, $\xi_{A_v}$ is independent of $\xi_{V \setminus B_v}$, $\xi_{B_v}$ is independent of $\xi_{V \setminus C_v}$ and $\mathrm{card}(C_v) \leq \kappa$. Then for all $x \in \mathbb{R}$,*

$$\sup_{x \in \mathbb{R}} |P[\xi_V \leq x] - \Phi(x)| \leq 75\kappa^2 \sum_{v \in V} E[|\xi_v|^3].$$

PROOF OF THEOREM 7. The argument is related to that used in [1] and elsewhere for various geometrical central limit theorems, but here has no geometric content.

Set $m = \mathrm{card}(V)$, and set

$$r := \delta \log m \qquad \text{with } \delta := \frac{\varepsilon}{12 \log D}.$$

For $v \in V$, set $\mathcal{N}_r(v) := \mathcal{N}_r(v, W)$ [as defined at (4)], and set

$$S_t := S_t(V, W) = \sum_{v \in V} X_t^W(v), \qquad \widetilde{S}_t := \sum_{v \in V} X_t^{\mathcal{N}_r(v)}(v),$$

$$\sigma := \sqrt{\mathrm{Var}(S_t)}, \qquad \tilde{\sigma} := \sqrt{\mathrm{Var}(\widetilde{S}_t)}.$$

We start by observing that provided we take $C(D, \varepsilon) \geq 1$, the estimate (13) will trivially hold if $\sigma < 1$. Therefore, in the sequel we may assume with no loss of generality that $\sigma \geq 1$.

For any two subsets $U, U'$ of $V$, it is the case that $\sum_{v \in U} X_t^{\mathcal{N}_r(v)}(v)$ and $\sum_{v \in U'} X_t^{\mathcal{N}_r(v)}(v)$ are independent if the graph distance between $U$ and $U'$ exceeds $2r$, because by definition the value of $X_t^{\mathcal{N}_r(v)}(v)$ is determined by $\{\tau_w : w \in \mathcal{N}_r(v)\}$ so that $\sum_{v \in U} X_t^{\mathcal{N}_r(v)}(v)$ is determined by the variables $\tau_w$, $w$ within distance $r$ of $U$, and likewise for $U'$.

For $v \in V$, set $A_v := \mathcal{N}_{2r}(v)$, $B_v := \mathcal{N}_{4r}(v)$ and $C_v := \mathcal{N}_{6r}(v)$. Also, set

$$\xi_v := \tilde{\sigma}^{-1}(X_t^{\mathcal{N}_r(v)}(v) - EX_t^{\mathcal{N}_r(v)}(v)).$$

Then the conditions for Lemma 14 hold, with

$$\kappa = \max_{v \in V}(\text{card}(C_v)) \leq (1 + D + D^2 + \cdots + D^{6r}) \leq D^{6r+1}$$

and $E[|\xi_v|^3] \leq \tilde{\sigma}^{-3}$. Applying Lemma 14 gives us

(69)
$$\sup_{x \in \mathbb{R}} |P[(\tilde{\sigma}^{-1}(\tilde{S}_t - E\tilde{S}_t)) \leq x] - \Phi(x)| \leq 75 D^{12r+2} \tilde{\sigma}^{-3} m$$
$$\leq 75 D^2 m^{1+\varepsilon} \tilde{\sigma}^{-3}.$$

We need to compare $S_t$ to $\tilde{S}_t$ and $\sigma$ to $\tilde{\sigma}$. By the proof of Proposition 1, for any $v \in V$ we have

$$P[X_t^W(v) \neq X_t^{\mathcal{N}_r(v)}(v)] \leq D^r/r!$$

(this holds for any of blocking RSA, dimer RSA or the annihilation process). Hence

(70)
$$P[S_t \neq \tilde{S}_t] \leq E[|S_t - \tilde{S}_t|] \leq mD^r/r! \leq K_1^r/r!$$

for $K_1 := e^{1/\delta} D$. Also, if we set $Y_t(v) = (X_t^{\mathcal{N}_r(v)}(v) - X_t^W(v))$ we have $|Y_t(v)| \leq 1$ with $P[Y_t(v) \neq 0] \leq D^r/r!$ so that

$$|\text{Cov}(S_t, \tilde{S}_t - S_t)| = \left|\sum_{v \in V} \sum_{w \in V} \text{Cov}(X_t^W(v), Y_t(w))\right| \leq 2m^2 D^r/r! \leq 2K_2^r/r!,$$

where we set $K_2 := e^{2/\delta} D$, and similarly

$$|\text{Var}(\tilde{S}_t - S_t)| = \left|\sum_{v \in V} \sum_{w \in V} \text{Cov}(Y_t(v), Y_t(w))\right| \leq 2m^2 D^r/r! \leq 2K_2^r/r!.$$

Hence,

$$|\tilde{\sigma}^2 - \sigma^2| = |2\text{Cov}(S_t, \tilde{S}_t - S_t) + \text{Var}(\tilde{S}_t - S_t)| \leq 6K_2^r/r!,$$

so that

$$\left|\left(\frac{\tilde{\sigma}}{\sigma}\right)^2 - 1\right| \leq 6\sigma^{-2} K_2^r/r!.$$

By the mean value theorem, for $\alpha \in \mathbb{R}$ there exists $C = C(\alpha)$ such that

(71)
$$\left|\left(\frac{\tilde{\sigma}}{\sigma}\right)^\alpha - 1\right| \leq \max(|(1 + 6\sigma^{-2} K_2^r/r!)^{\alpha/2} - 1|, |(1 - 6\sigma^{-2} K_2^r/r!)^{\alpha/2} - 1|)$$
$$\leq C\sigma^{-2} K_2^r/r!.$$

Now consider arbitrary $x \in \mathbb{R}$ and set $y := x\sigma/\tilde{\sigma}$. Then by the mean value theorem, with $\phi$ denoting the standard normal probability density function, we have

$$\begin{aligned}|\Phi(y) - \Phi(x)| &\leq |y - x|\phi(\min(|x|, |y|)) \\ &= (\max(|x|, |y|) - \min(|x|, |y|))\phi(\min(|x|, |y|)) \\ &= \max\left(\frac{\sigma}{\tilde{\sigma}} - 1, \frac{\tilde{\sigma}}{\sigma} - 1\right)\min(|x|, |y|)\phi(\min(|x|, |y|)) \\ &\leq C\sigma^{-2}K_2^r/r!,\end{aligned}$$

where we have used (71) and the fact that $\max_{u \geq 0}\{u\phi(u)\} \leq 1$. Using (69) and (71), we then obtain

$$\begin{aligned}(72) \quad &|P[\sigma^{-1}(\widetilde{S}_t - E\widetilde{S}_t) \leq x] - \Phi(x)| \\ &\leq |P[\tilde{\sigma}^{-1}(\widetilde{S}_t - E\widetilde{S}_t) \leq y] - \Phi(y)| + |\Phi(y) - \Phi(x)| \\ &\leq 75D^2 m^{1+\varepsilon}\sigma^{-3}(1 + C\sigma^{-2}K_2^r/r!) + C\sigma^{-2}K_2^r/r! \\ &\leq C'm^{1+\varepsilon}\sigma^{-3},\end{aligned}$$

where in the last line we have used the assumption that $\sigma \geq 1$ and the (easily proved) fact that $\sigma \leq m$.

Next, for $x \in \mathbb{R}$ we use (70) and (72) and the fact that $\sigma \leq m$ to obtain

$$\begin{aligned}&|P[\sigma^{-1}(\widetilde{S}_t - ES_t) \leq x] - \Phi(x)| \\ &\leq |\Phi(x + \sigma^{-1}(ES_t - E\widetilde{S}_t)) - \Phi(x)| \\ &\quad + |P[\sigma^{-1}(\widetilde{S}_t - E\widetilde{S}_t) \leq x + \sigma^{-1}(ES_t - E\widetilde{S}_t)] \\ &\qquad - \Phi(x + \sigma^{-1}(ES_t - E\widetilde{S}_t))| \\ &\leq (\sigma^{-1}K_1^r/r!) + C'm^{1+\varepsilon}\sigma^{-3} = m\sigma^{-3}(C'm^{\varepsilon} + (K_1^r/r!)\sigma^2/m) \\ &\leq m\sigma^{-3}(C'm^{\varepsilon} + mK_1^r/r!) = m\sigma^{-3}(C'm^{\varepsilon} + K_2^r/r!),\end{aligned}$$

since $K_2 = K_1 e^{1/\delta}$. Hence, (72) remains valid with $E\widetilde{S}_t$ in the left-hand side replaced by $ES_t$ (and also with a change of constant). Finally, by (70) and the fact that $\sigma \leq m$,

$$\begin{aligned}&m^{-1}\sigma^3|P[\sigma^{-1}(S_t - ES_t) \leq x] - P[\sigma^{-1}(\widetilde{S}_t - ES_t) \leq x]| \\ &\leq m^{-1}\sigma^3 P[S_t \neq \widetilde{S}_t] \leq C'(K_3^r/r!)\end{aligned}$$

(with $K_3 := e^{3/\delta}D$) which is also bounded. Hence, having replaced $E\widetilde{S}_t$ by $ES_t$ in (72), we can also replace $\widetilde{S}_t$ by $S_t$, completing the proof. $\square$

PROOF OF THEOREM 8. Let $0 < t < 1$, let $W$ be the vertex set of a graph with degrees bounded by $D$ and let $V \subseteq W$ be finite; set $m = \text{card}(V)$.

First consider blocking RSA. Choose an independent set $I$ in $V$ with at least $m/(D+1)$ elements; this is possible by a greedy algorithm.

For $v \in I$, let $Y_v$ denote the indicator of the event that $\tau_w > t$ for $w \in \mathcal{N}_1(v, W) \setminus \{v\}$. Then

$$E \sum_{v \in I} Y_v \geq \left(\frac{m}{D+1}\right)(1-t)^D. \tag{73}$$

Let $\mathcal{F}$ be the $\sigma$-field generated by the values of $Y_v, v \in I$ along with the values of

$$\tau_w, w \in W \setminus \left( \bigcup_{\{v \in I : Y_v = 1\}} (\mathcal{N}_1(v)) \right).$$

Then, setting $S_t := S_t(V, W)$, we have

$$\begin{aligned} \operatorname{Var}(S_t) &= \operatorname{Var}(E[S_t | \mathcal{F}]) + E[\operatorname{Var}(S_t | \mathcal{F})] \\ &\geq E[\operatorname{Var}(S_t | \mathcal{F})]. \end{aligned} \tag{74}$$

If we are given the values of $Y_v, v \in I$, along with the values of $\tau_w$ for $w \in W \setminus (\bigcup_{\{v \in I : Y_v = 1\}} (\mathcal{N}_1(v)))$, then the remaining variability of $S_t$ comes only from the values of $X_t(v)$ for $v \in I$ with $Y_v = 1$. For such $v$, site $v$ contributes 1 to $S_t$ if $\tau_v \leq t$ and contributes 0 to $S_t$ if $\tau_v > t$. Hence,

$$\operatorname{Var}(S_t | \mathcal{F}) = \operatorname{Var}\left( \operatorname{Bin}\left( \sum_{v \in I} Y_v, t \right) \right) = t(1-t) \sum_{v \in I} Y_v$$

so that by (73) and (74),

$$\operatorname{Var}(S_t) \geq t(1-t) E \sum_{v \in I} Y_v \geq t(1-t) m (1-t)^D / (D+1),$$

which completes the proof of (14) for blocking RSA.

The argument is similar for dimer RSA. Choose an independent set $I'$ of *edges*, each incident to at least one element of $V$, with at least $m/(2D)$ elements of $I'$, using a greedy algorithm and the assumption that there are no isolated vertices in $V$. For each edge $e \in I'$, let $Y_v$ be the indicator of the event that $\tau_f > t$ for all edges $f$ adjacent to $e$. The rest of the argument proceeds much as before.

A similar argument applies for the annihilation process. The factor of $(1/2)$ is included in the right-hand side of (15) because edges in $I'$ might have only one endpoint in $V$. We do not give details for this case. $\square$

PROOF OF THEOREM 9. Suppose $V \subseteq W$, where $V$ is finite and $W$ is the finite or countable vertex set of a graph with all degrees at most $D$. Then $\operatorname{card}(\mathcal{N}_5(v, W)) \leq D^6$ for all $v \in V$. Recall that $W^+$ denotes the set of sites in $W$ with positive entropy. Using a greedy algorithm and the bound on degrees, choose a collection $J$ of sites in $V \cap W^+$, all lying a graph distance at least 6 from one other, and satisfying

$$\text{card}(J) > \frac{\text{card}(V \cap W^+)}{D^6}. \tag{75}$$

For $v \in J$, let $Z_v$ be the indicator of the event that for all sites $w$ in the associated blocking set $B_v$ (see the definition of positive entropy in Section 2.7), we have $\tau_w \leq 1/2$, and $\tau_x > 1/2$ for $x \in (\mathcal{N}_1(w, W)) \setminus \{w\}$. Then

$$P[Z_v = 1] \geq (\tfrac{1}{2})^{D^3(1+D)}. \tag{76}$$

Let $\mathcal{F}'$ be the $\sigma$-field generated by the values of $Z_v, v \in J$, along with the values of

$$\tau_w, w \in W \setminus \left( \bigcup_{\{v \in J \,:\, Z_v = 1\}} (\mathcal{N}_1(v) \cup \{v\}) \right).$$

Then, as at (74), setting $S_1 := S_1(V, W)$ we have

$$\text{Var}(S_1) \geq E[\text{Var}(S_1 | \mathcal{F}')]. \tag{77}$$

If we are given the values of $Z_v, v \in J$, along with the values of $\tau_w$ for $w \in W \setminus (\bigcup_{\{v \in J \,:\, Z_v = 1\}} (\mathcal{N}_1(v) \cup \{v\}))$, then the remaining variability of $S_1$ comes only from the values of $\{X_w, w \in \mathcal{N}_1(v, W)\}$ with $v \in J$ and $Z_v = 1$. For such $v$, pick $u, u'$ adjacent to $v$ but not to each other. The collection of sites $\mathcal{N}_1(v)$ contributes 1 to $S_1$ if $\tau_v \leq \min_{w \in \mathcal{N}_1(v) \setminus \{v\}} \tau_w$ [an event of probability at least $1/(D+1)$], and contributes at least 2 to $S_1$ if $\max(\tau_u, \tau_{u'}) < \min\{\tau_y : y \in \mathcal{N}_1(v) \setminus \{u, u'\}\}$ [an event of probability at least $2(D+1)^{-2}$]. For any random variable $X$ we have $\text{Var}(X) \geq (1/4) \min(P[X \geq 2], P[X = 1])$. Hence,

$$\text{Var}(S_1 | \mathcal{F}) \geq \tfrac{1}{2}(D+1)^{-2} \sum_{v \in J} Z_v$$

so that by (75) and (76),

$$\text{Var}(S_1) \geq \frac{1}{2}(D+1)^{-2} \sum_{v \in J} E[Z_v] \geq \frac{(1/2)^{D^3(1+D)+1} \text{card}(V \cap W^+)}{D^6(D+1)^2},$$

which completes the proof of (16) for blocking RSA. The proof of (17) for the annihilation process is similar. $\square$

PROOF OF THEOREM 11. By Theorem 7, it suffices for us to prove the convergence of variances (19) and (20) for blocking RSA, along with the corresponding convergence of variances for dimer RSA and for the annihilation process.

DEPOSITION PROCESSES ON GRAPHS 35We first show the equivalence of (19) and (20). Let $\Delta_i^n := X_t^{L_n}(i) - X_t^{\mathbb{Z}}(i)$. Then $|\Delta_i^n| \leq 1$ and by the argument in the proof of Proposition 1, for any of blocking or dimer RSA or the annihilation process we have

$$P[\Delta_i^n \neq 0] \leq \frac{1}{i!} + \frac{1}{(n+1-i)!}.$$

Since $|\Delta_i^n| \leq 1$, we have

$$|\operatorname{Cov}(\Delta_i^n, \Delta_j^n)| \leq E|\Delta_i^n \Delta_j^n| + E|\Delta_i^n|E|\Delta_j^n| \leq 2\min(P[\Delta_i^n \neq 0], P[\Delta_j^n \neq 0]).$$

Hence $\operatorname{Var} \sum_{i=1}^n \Delta_i^n = O(1)$ since

$$\operatorname{Var} \sum_{i=1}^n \Delta_i^n = \sum_{i=1}^n \sum_{j=1}^n \operatorname{Cov}(\Delta_i^n, \Delta_j^n)$$

$$\leq 2 \sum_{i=1}^n \sum_{j=1}^n \min\left(\frac{1}{i!} + \frac{1}{(n+1-i)!}, \frac{1}{j!} + \frac{1}{(n+1-j)!}\right)$$

$$\leq 8 \sum_{i=1}^{\lceil n/2 \rceil} \sum_{j=1}^{\lceil n/2 \rceil} \min\left(\frac{2}{i!}, \frac{2}{j!}\right)$$

$$\leq 32 \sum_{i=1}^{\infty} \frac{1}{(i-1)!}$$

$$= 32 e^{-1}.$$

Hence, if (19) holds, then by the Cauchy–Schwarz inequality

$$\operatorname{Cov}\left(\sum_{i=1}^n \Delta_i^n, S_t(L_n, \mathbb{Z})\right) = O(n^{1/2})$$

and since $S_t(L_n, L_n) = S_t(L_n, \mathbb{Z}) + \sum_{i=1}^n \Delta_i^n$, we obtain

$$\operatorname{Var} S_t(L_n, L_n) - \operatorname{Var} S_t(L_n, \mathbb{Z})$$

(78)
$$= \operatorname{Var} \sum_{i=1}^n \Delta_i^n + 2\operatorname{Cov}\left(\sum_{i=1}^n \Delta_i^n, S_t(L_n, \mathbb{Z})\right)$$

$$= O(n^{1/2}),$$

so that (19), if true, implies (20). Similarly (20), if true, implies (19).

Define $\operatorname{Cov}(m, t) := \operatorname{Cov}(m, t, T_1)$ as given by (27), with $\operatorname{Cov}(0, t) = \operatorname{Var} X_t^{\mathbb{Z}}(0)$. Then by translation invariance,

$$n^{-1} \operatorname{Var}(S_t(L_n, \mathbb{Z})) = n^{-1} \sum_{i=1}^n \sum_{j=1}^n \operatorname{Cov}(|i-j|, t)$$

$$
(79) \qquad = \operatorname{Cov}(0,t) + 2n^{-1} \sum_{m=1}^{n-1} (n-m) \operatorname{Cov}(m,t)
$$

$$
\to \operatorname{Cov}(0,t) + 2 \sum_{m=1}^{\infty} \operatorname{Cov}(m,t),
$$

where the last line follows from the dominated convergence theorem.

In the case of blocking RSA, using Theorem 13 along with the case $k=1$ of (29), and collecting terms, we obtain

$$
\sum_{r=1}^{\infty} \operatorname{Cov}(r,t) = -\frac{e^{-2t}}{2} \sum_{r=1}^{\infty} \sum_{n=0}^{\infty} \frac{(-2t)^{r+1+2n}}{(r+1+2n)!}
$$

$$
(80) \qquad = -\frac{e^{-2t}}{2} \left( \sum_{\text{even } m>0} \left(\frac{m}{2}\right) \frac{(-2t)^m}{m!} + \sum_{\text{odd } m>0} \left(\frac{m-1}{2}\right) \frac{(-2t)^m}{m!} \right)
$$

$$
= -\frac{e^{-2t}}{4}(-2te^{-2t} - \sinh(-2t)).
$$

Also, by Theorem 12, $\operatorname{Cov}(0,t)$ is the variance of a Bernoulli variable with parameter $(1+e^{-2t})/2$, so that

$$
\operatorname{Cov}(0,t) = \frac{1}{4}(1+e^{-2t})(1-e^{-2t}) = \frac{1-e^{-4t}}{4} = e^{-2t}\sinh(2t)/2,
$$

and combining this with (80), we see that the limit in (79) equals $te^{-4t}$, which gives us (19) as required for monomer RSA.

In the case of dimer RSA, by using (40) and (46), and some routine algebra which we omit, we obtain

$$
\operatorname{Cov}(0,t) + 2 \sum_{m=1}^{\infty} \operatorname{Cov}(m,t) = 4te^{-4t}
$$

so that by (79), $\sigma_{\mathrm{DI}}^2 = 4te^{-4t}$ as asserted. In fact, this can also be derived from the blocking RSA result by a duality argument; dimer RSA on a row of $n$ vertices is equivalent to blocking RSA on a row of $n-1$ vertices, where each site for the blocking process corresponds to a bond (edge) for the dimer process, and hence an occupied site in the blocking process counts as a pair of occupied sites in the dimer process, so that $S_t(L_n, L_n)$ for the dimer process has the same distribution as $2S_t(L_{n-1}, L_{n-1})$ for the blocking process, and hence $\sigma_{\mathrm{DI}}^2(t) = 4\sigma_{\mathrm{BL}}^2(t)$.

For the annihilation process, using (48) and (50) along with (79) we have

$$
\sigma_{\mathrm{AN}}^2(t) = \operatorname{Cov}(0,t) + 2 \sum_{m=1}^{\infty} \operatorname{Cov}(m,t) = e^{-t}(1-e^{-t}) - 2e^{-t} \sum_{m=1}^{\infty} \sum_{j=m+1}^{\infty} \frac{(-t)^j}{j!}
$$

$$= e^{-t} - e^{-2t} - 2e^{-t}\bigg(\sum_{m=1}^{\infty}\sum_{j=m}^{\infty}\frac{(-t)^j}{j!} - \sum_{m=1}^{\infty}\frac{(-t)^m}{m!}\bigg)$$

$$= e^{-t} - e^{-2t} + 2e^{-t}\bigg(e^{-t} - 1 - \sum_{j=1}^{\infty}\frac{(-t)^j}{(j-1)!}\bigg)$$

$$= e^{-2t} - e^{-t} + 2te^{-2t},$$

as asserted. $\square$

**Acknowledgment.** We thank the referee and an Associate Editor for some helpful comments.

## REFERENCES

[1] AVRAM, F. and BERTSIMAS, D. (1993). On central limit theorems in geometric probability. *Ann. Appl. Probab.* **3** 1033–1046. MR1241033
[2] BRIGHTWELL, G. R., HÄGGSTRÖM, O. and WINKLER, P. (1999). Nonmonotonic behavior in hard-core and Widom–Rowlinson models. *J. Statist. Phys.* **94** 415–435. MR1675359
[3] CASER, S. and HILHORST, H. J. (1994). Exact bounds from a new series-expansion method for random sequential adsorption. *J. Phys. A* **27** 7969–7979. MR1323573
[4] CHEN, L. H. Y. and SHAO, Q.-M. (2004). Normal approximation under local dependence. *Ann. Probab.* **32** 1985–2028. MR2073183
[5] EVANS, J. W. (1984). Exactly solvable irreversible processes on Bethe lattices. *J. Math. Phys.* **25** 2527–2532. MR751543
[6] EVANS, J. W. (1989). Comment on "Kinetics of random sequential adsorption." *Phys. Rev. Lett.* **62** 2642.
[7] EVANS, J. W. (1993). Random and cooperative sequential adsorption. *Rev. Modern Phys.* **65** 1281–1329.
[8] FAN, Y. and PERCUS, J. K. (1991). Asymptotic coverage in random sequential adsorption on a lattice. *Phys. Rev. A* **44** 5099–5103. MR1133846
[9] GALVIN, D. and KAHN, J. (2004). On phase transition in the hard-core model on $\mathbb{Z}^d$. *Combin. Probab. Comput.* **13** 137–164. MR2047233
[10] GRIMMETT, G. (1999). *Percolation*, 2nd ed. Springer, New York. MR1707339
[11] LIGGETT, T. (1985). *Interacting Particle Systems.* Springer, Berlin. MR776231
[12] LYONS, R. and PERES, Y. (2004). *Probability on Trees and Networks.* Book in preparation. Available at http://mypage.iu.edu/~rdlyons/.
[13] MEAKIN, P., CARDY, J. L., LOH, E., JR. and SCALAPINO, D. J. (1986). Maximal coverage in random sequential adsorption. *J. Chem. Phys.* **86** 2380–2382.
[14] O'HELY, M. and SUDBURY, A. (2001). The annihilating process. *J. Appl. Probab.* **38** 223–231. MR1816125
[15] PAGE, E. S. (1959). The distribution of vacancies on a line. *J. Roy. Statist. Soc. Ser. B* **21** 364–374. MR119218
[16] PEDERSEN, F. B. and HEMMER, P. (1993). Time evolution of correlations in a random sequential adsorption process. *J. Chem. Phys.* **98** 2279–2282.
[17] PEMANTLE, R. and STACEY, A. (2001). Branching random walk and contact process on Galton–Watson and non-homogeneous trees. *Ann. Probab.* **29** 1563–1590. MR1880232

[18] PENROSE, M. D. (2001). Random parking, sequential adsorption, and the jamming limit. *Comm. Math. Phys.* **218** 153–176. MR1824203
[19] PENROSE, M. D. (2002). Limit theorems for monotonic particle systems and sequential deposition. *Stochastic Process. Appl.* **98** 175–197. MR1887532
[20] PENROSE, M. D. and YUKICH, J. E. (2002). Limit theorems for random sequential packing and deposition. *Ann. Appl. Probab.* **12** 272–301. MR1890065
[21] SUDBURY, A. (2002). Inclusion–exclusion methods for treating annihilating and deposition processes. *J. Appl. Probab.* **39** 466–478. MR1928883
[22] VAN DEN BERG, J. and STEIF, J. E. (1994). Percolation and the hard-core lattice gas model. *Stochastic Process. Appl.* **49** 179–197. MR1260188
[23] WIDOM, B. (1966). Random sequential addition of hard spheres to a volume. *J. Chem. Phys.* **44** 3888–3894.
[24] WIDOM, B. (1973). Random sequential filling of intervals on a line. *J. Chem. Phys.* **58** 4043–4044.

DEPARTMENT OF MATHEMATICAL SCIENCES
UNIVERSITY OF BATH
BATH BA2 7AY
UNITED KINGDOM
E-MAIL: m.d.penrose@maths.bath.ac.uk

SCHOOL OF MATHEMATICAL SCIENCES
MONASH UNIVERSITY
VICTORIA 3800
AUSTRALIA
E-MAIL: aidan.sudbury@sci.monash.edu.au